\documentclass{article}

\usepackage{amsmath,amssymb,amsthm,enumerate,mathrsfs}
\usepackage[numbers]{natbib}
\usepackage{booktabs}
\usepackage{adjustbox}
\usepackage{subfigure,multirow}
\usepackage{appendix}
\usepackage{authblk}
\usepackage{geometry}
\usepackage[ruled,linesnumbered]{algorithm2e}
\usepackage{algpseudocode}%
\usepackage{longtable}
\usepackage{caption}
\usepackage{lipsum}
\usepackage{titlesec}

\newcommand{\x}{\mathbf{x}}
\newcommand{\yy}{\mathbf{y}}
\newcommand{\sss}{\mathbf{s}}
\newcommand{\g}{\mathbf{g}}

\newcommand{\eqf}{equilibrium }
\newcommand{\xo}{\mathbf{0}}

\newcommand{\T}{\text{T}}
\newcommand{\BB}{Barzilai-Borwein }

\theoremstyle{plain}
\newtheorem{theorem}{Theorem}
\newtheorem{lemma}[theorem]{Lemma}

\theoremstyle{remark}
\newtheorem{remark}{Remark}

\newtheorem{property}{Property}

\title{\bf On Convergence and Stability of Two Extended BB-like Step Sizes\thanks{The work is supported by the National Natural Science Foundation of China (Project No. 12371099).}}

\author{Xin Xu\thanks{Corresponding author: xustonexin@gmail.com}}
\affil{School of Mathematics, \\Southwestern University of Finance and Economics, Chengdu, China}
\date{}

\begin{document}
\maketitle

\begin{abstract}
	The \BB(BB) step sizes have a profound impact on gradient descent methods. In this work, we propose two new gradient step sizes: one longer than the original long BB step size, and the other shorter than the original short BB step size. This extends the bounds of the original BB step sizes. For strictly convex quadratic optimization problems, we employ the dynamics of difference equations to prove that these two new methods achieve R-linear convergence. Regarding stability, we surprisingly find that under certain conditions, the gradient descent method corresponding to the longer step size is stable, whereas the shorter step size consistently leads to instability. Numerical results validate these stability theories. Here, stability refers to whether the gradient norm decreases monotonically.
	
	
\end{abstract}

\section{Introduction}\label{sec:introduction}
The gradient descent method for unconstrained minimization of $f:\mathbb{R}^{n}\rightarrow\mathbb{R}$ takes the form
\begin{equation}\label{equ:BB}
\x_{k+1}=\x_{k}-\alpha_{k}^{-1}\g_{k},
\end{equation}
where $\g_{k}=\nabla f(\x_{k})$ denotes the gradient at $\x_{k}$, and $\alpha_{k}^{-1}>0$ is the step size. 

The classical \BB (BB) method \cite{Barzilai1988TwoPointStep} determines two scalars, denoted 
\begin{equation}\label{BB steps}
\alpha_{k}^{BB1}=\frac{\sss_{k-1}^{\T}\yy_{k-1}}{\sss_{k-1}^{\T}\sss_{k-1}}\quad\text{and}\quad \alpha_{k}^{BB2}=\frac{\yy_{k-1}^{\T}\yy_{k-1}}{\sss_{k-1}^{\T}\yy_{k-1}},
\end{equation}
which are solutions to the following least-squares problems
\begin{equation}\label{LS}
\min_{\alpha\in\mathbb{R}} \|\alpha \sss_{k-1} - \yy_{k-1}\|_{2}^{2}\quad\text{and}\quad \min_{\alpha\in \mathbb{R}} \|\sss_{k-1} -  {\alpha}^{-1}\yy_{k-1}\|_{2}^{2},
\end{equation} 
where $\sss_{k-1}=\x_{k}-\x_{k-1}$ and $\yy_{k-1}=\g_{k}-\g_{k-1}$. When $\sss_{k-1}^{\T}\yy_{k-1}>0$, these satisfy $\alpha_{k}^{BB1}\le\alpha_{k}^{BB2}$. In the  case of convex quadratic minimization problems
\begin{equation}\label{equ:qua}
\min_{\x\in\mathbb{R}^{n}} f(\x)=\frac{1}{2}\x^{\T}A\x-\mathbf{b}^{\T}\x,
\end{equation}
where $A\in\mathbb{R}^{n\times n}$ is a symmetric positive definite (SPD) matrix and $\mathbf{b}\in\mathbb{R}^{n}$, the BB scalars are
\begin{equation}\label{BB stepsA}
\alpha_{k}^{BB1}=\frac{\sss_{k-1}^{\T}A\sss_{k-1}}{\sss_{k-1}^{\T}\sss_{k-1}}\quad\text{and}\quad \alpha_{k}^{BB2}=\frac{\sss_{k-1}^{\T}A^2\sss_{k-1}}{\sss_{k-1}^{\T}A\sss_{k-1}}.
\end{equation}
The BB method approximates the secant equation by the scalar matrices $\alpha_{k}^{BB1}I$ or $\alpha_{k}^{BB2}I$, where $I$ is the identity matrix, thereby classifying it as a quasi-Newton method. Owing to its simplicity and effectiveness, the BB method has been the subject of extensive analysis and numerous enhancements within the optimization literature, such as \cite{Raydan1993BarzilaiBorweinchoice,Dai2002Rlinearconvergence,Fletcher2005BarzilaiBorweinMethod,Zhou2006GradientMethodsAdaptive,Frassoldati2008Newadaptivestepsize,Bonettini2009scaledgradientprojection,Dai2013NewAnalysisBarzilai,An2025RegularizedBarzilaiBorwein}.

Recently, \cite{Xu2025ParameterizedBarzilaiBorwein} introduced a Parametrized BB (PBB) method. This approach is based on an interpolated least-squares model
\begin{equation}\label{VLS}
\min_{\alpha\in\mathbb{R}} L(\alpha):=\|\alpha^{m_k} \sss_{k-1} -\alpha^{m_k-1} \yy_{k-1}\|_{2}^{2},
\end{equation}
yielding the scalar
\begin{equation}\label{VBB}
\alpha_{k}^{PBB}=\frac{(2m_k-1)\sss_{k-1}^{\T}\yy_{k-1}+\sqrt{\big((2m_k-1)\sss_{k-1}^{\T}\yy_{k-1}\big)^{2}-4m_k(m_k-1)\sss_{k-1}^{\T}\sss_{k-1}\yy_{k-1}^{\T}\yy_{k-1}}}{2m_k\sss_{k-1}^{\T}\sss_{k-1}},
\end{equation} 
where $m_k\in(0, 1]$. It follows from Theorem 1 in \cite{Xu2025ParameterizedBarzilaiBorwein} that $\alpha_{k}^{PBB}\in[\alpha_{k}^{BB1},\alpha_{k}^{BB2})$, and that $\alpha_{k}^{PBB}$ is monotonically decreasing in $m_{k}$. By selecting an interpolation parameter $m_{k}$, $\alpha_{k}^{PBB}$ interpolates between $\alpha_{k}^{BB1}$ and $\alpha_{k}^{BB2}$ to adapt to problem characteristics. The PBB formulation provides a concise unification of the original BB step sizes, with $\alpha_{k}^{BB1}$ and $\alpha_{k}^{BB2}$ representing its extreme cases. A limitation of PBB, however, is that it does not extend beyond the range bounded by $\alpha_{k}^{BB1}$ and $\alpha_{k}^{BB2}$, resulting in $\alpha_{k}^{PBB}$ being a quasi-convex combination of them. 

In this work, building upon the interpolated least-squares model \eqref{VLS}, we propose two new scalars: one less than $\alpha_{k}^{BB1}$ and another greater than $\alpha_{k}^{BB2}$. The resulting step sizes thus extend beyond the conventional BB bounds, offering enhanced flexibility for gradient descent algorithms. More significantly, our theoretical analysis reveals that when certain conditions are met, gradient descent employing the longer step size exhibits monotonic gradient norm decrease. Conversely, the method employing the shorter step size is inherently non-monotonic. These insights facilitate the design of more efficient gradient descent schemes. The derivation of the new step sizes is presented in Section 2. For the two-dimensional quadratic minimization problem \eqref{equ:qua}, we analyze the convergence and stability properties of the corresponding methods from the perspective of difference equation dynamics in Section 3. In Section 4, we propose a truncation mechanism that restricts the proposed scalars within the spectral radius of the Hessian matrix. Unless otherwise stated, throughout this paper $\|\cdot\|_{2}$ denotes the Euclidean $2$-norm of vectors, $\mathbf{1}\in\mathbb{R}^{n}$ is the all-ones vector, $\mathbf{0}\in\mathbb{R}^{n}$ is the all-zeros vector, and $\kappa(A)=\frac{\lambda_{1}}{\lambda_{n}}$ denotes the spectral condition number of an $n\times n$  matrix $A$, where $\lambda_{1}$ and $\lambda_{n}$ are its largest and smallest eigenvalues, respectively. The spectral radius of $A$ is defined as
\begin{equation*}
\rho(A)=\max\{|\lambda|:\lambda\in\sigma(A)\},
\end{equation*}
where $\sigma(A)$ denotes the spectrum of $A$ ( the set of all eigenvalues).

\section{Derivation of new step sizes}\label{sec:twostep}
In the PBB method, $m_{k}\in[0,1]$ yields  $\alpha_{k}^{PBB}\in[\alpha_{k}^{BB1},\alpha_{k}^{BB2}]$. As discussed in Section \ref{sec:introduction}, our objective is to extend the bounds of the original BB step sizes to provide more effective choice for gradient descent. To this end, an intuitive approach is to expand the range of $m_{k}$.

Considering the problem \eqref{VLS}, its first-order optimality condition reduces, after simplification, to the quadratic equation
\begin{equation}\label{equ:ppeizhi}
\phi(\alpha):=m_k\sss_{k-1}^{\T}\sss_{k-1}\alpha^{2}-(2m_k-1)\sss_{k-1}^{\T}\yy_{k-1}\alpha+(m_k-1)\yy_{k-1}^{\T}\yy_{k-1}=0,
\end{equation}
in the variable $\alpha_{k}$ when $m_{k}\neq 0$. The discriminant of the quadratic function $\phi(\alpha)$ is given by   $$\Delta=(2m_{k}-1)^2(\sss_{k-1}^{\T}\yy_{k-1})^2-4m_{k}(m_{k}-1)\sss_{k-1}^{\T}\sss_{k-1}\yy_{k-1}^{\T}\yy_{k-1}.$$ 
Define the angle between vectors  
$$\theta_{k}=\angle(\sss_{k-1}, \yy_{k-1}),$$
which yields the trigonometric relation 
\begin{equation}\label{equ:cosk}
\cos^{2}\theta_{k}=\frac{(\sss_{k-1}^{\T}\yy_{k-1})^2}{\|\sss_{k-1}\|_{2}^{2}\|\yy_{k-1}\|_{2}^{2}}.
\end{equation} 
Setting 
\begin{equation*}\label{equ:n}
n_k=\frac{\cos^2\theta_{k}}{1-\cos^2\theta_{k}},
\end{equation*}
if either   
\begin{equation*}
m_{k1}=\frac{1+\sqrt{1+n_k}}{2}>1,
\end{equation*}
or
\begin{equation*}
m_{k2}=\frac{1-\sqrt{1+n_k}}{2}<0,
\end{equation*}
then $\Delta=0$, implying that  \eqref{equ:ppeizhi} admits a unique solution. This solution corresponds to  
\begin{equation}\label{equ:BB1left}
\alpha_{k}^{L}=\frac{\sss_{k-1}^{\T}\yy_{k-1}}{\sss_{k-1}^{\T}\sss_{k-1}(1+\sin\theta_{k})}\le\alpha_{k}^{BB1},
\end{equation}
or
\begin{equation}\label{equ:BB1right}
\alpha_{k}^{R}=\frac{\sss_{k-1}^{\T}\yy_{k-1}}{\sss_{k-1}^{\T}\sss_{k-1}(1-\sin\theta_{k})}\ge\alpha_{k}^{BB2},
\end{equation}
respectively. Remarkably, these solutions satisfy the relation 
\begin{equation*}\label{equ:mean1}
\alpha_{k}^{L}\alpha_{k}^{R}=\alpha_{k}^{BB1}\alpha_{k}^{BB2}.
\end{equation*}

\begin{figure}[h!]
	\centering
	\subfigure{\includegraphics[width=0.9\linewidth]{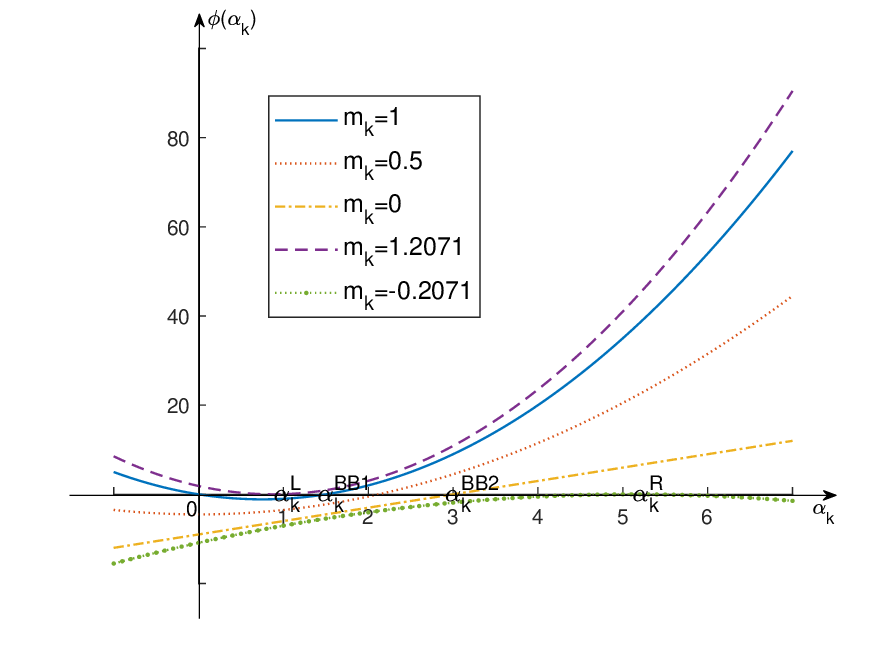}}
	\caption{\textit{The graphs of the $\phi(\alpha)$ quadratic functions  corresponding to different $m_{k}$.}}
	\label{fig:LRVBB}
\end{figure}

We now consider the quadratic function $\phi(\alpha)$ with the fixed parameter values:
\begin{equation*}
s_{k-1}^{\T}s_{k-1}=2,\quad s_{k-1}^{\T}y_{k-1}=3,\quad y_{k-1}^{\T}y_{k-1}=9.
\end{equation*}
Figure \ref{fig:LRVBB} displays $\phi(\alpha)$ for $m_{k}=m_{k1}$, $1$, $0.5$, $0$, $m_{k2}$, respectively. The positive roots along the $\alpha_{k}$-axis correspond to $\alpha_{k}^{L}$, $\alpha_{k}^{BB1}$, $\sqrt{\alpha_{k}^{BB1}\alpha_{k}^{BB2}}$,$\alpha_{k}^{BB2}$, and $\alpha_{k}^{R}$, respectively. As shown in Figure \ref{fig:LRVBB},  $\alpha_{k}^{L}$ lies to the left of $\alpha_{k}^{BB1}$, while $\alpha_{k}^{R}$ lies to the right of $\alpha_{k}^{BB2}$. We designate the gradient descent method with step size $(\alpha_{k}^{L})^{-1}$ as the LEFT method, and the method with step size $(\alpha_{k}^{R})^{-1}$ as the RIGHT method. This broadens the selectable range of step sizes for gradient descent.

\begin{figure}[h!]
	\centering
	\subfigure{\includegraphics[width=0.9\linewidth]{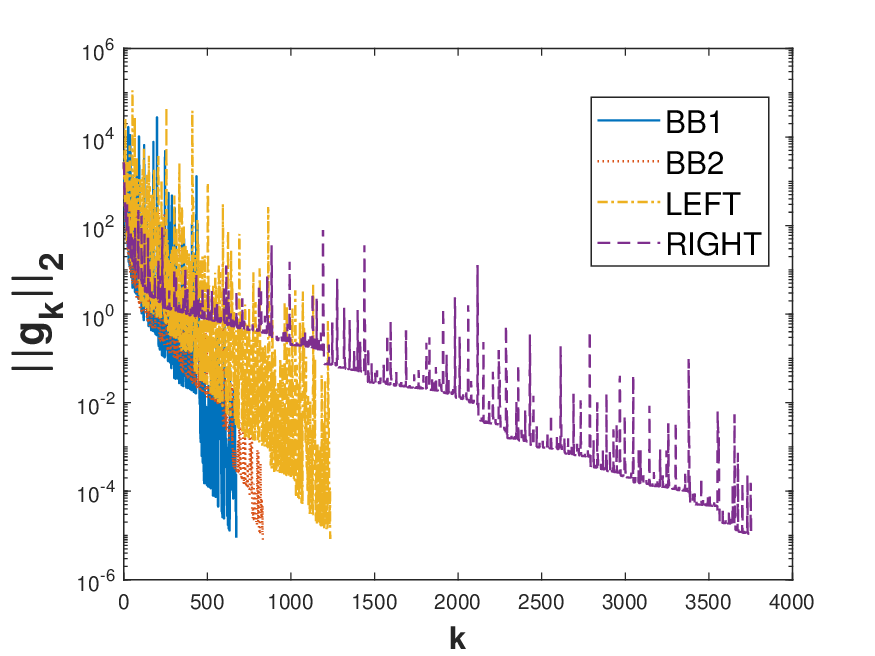}}
	\caption{\textit{The performance of the BB1, BB2, LEFT and RIGHT methods for problem \eqref{pro:qua}.}}
	\label{fig:LRVBBGnorm}
\end{figure}

Without loss of generality, we consider the following quadratic function 
\begin{equation}\label{pro:qua}
f(\mathbf{x})=\frac{1}{2}(\mathbf{x}-\mathbf{x}_{*})^{\T}A(\mathbf{x}-\mathbf{x}_{*}),
\end{equation}
where 
\begin{equation}\label{A}
A = \text{diag}\{\lambda_{1},\ldots,\lambda_{n}\}
\end{equation}
with eigenvalues ordered as
\begin{equation*}
\lambda=\lambda_{1}\ge\ldots\ge\lambda_{n-1}\ge\lambda_{n}=1.
\end{equation*}
Following the benchmark in \cite{DeAsmundis2014efficientgradientmethod}, we define the eigenvalues as
\begin{equation}\label{diagconditionn}
\lambda_{i}=10^{\frac{ncond}{n-1}(n-i)},\quad \text{for}\quad i=1,\ldots,n,
\end{equation}
where $ncond=\text{log}_{10}(\kappa(A))$. The solution is fixed at $\mathbf{x}_{*}=\mathbf{1}$ with initial guess $\mathbf{x}_{1}=\mathbf{0}$. We adopt dimension $n=10$, condition number $\kappa(A)=10^{4}$, and stopping criterion $\|\mathbf{g}_{k}\|_{2}\le 10^{-9}\Vert \mathbf{g}_{1}\Vert_{2}$. The initial step size is set to  $\frac{\mathbf{g}_{1}^{\T}\mathbf{g}_{1}}{\mathbf{g}_{1}^{\T}A\mathbf{g}_{1}}$. Figure \ref{fig:LRVBBGnorm} compares the performance of BB1, BB2, LEFT, and RIGHT methods on this quadratic problem using gradient norms $\|\g_{k}\|_{2}$.

Figure \ref{fig:LRVBBGnorm} demonstrates that both LEFT and RIGHT methods underperform the original BB1 and BB2 methods. The RIGHT method exhibits particularly poor efficiency, requiring more iterations to satisfy the termination criterion than BB2.

In Section \ref{sec:converge}, we conduct stability and convergence analyses of the LEFT and RIGHT methods to elucidate why these extended step sizes yield inferior performance compared to conventional BB step sizes. Building on these insights, we develop modified versions of the proposed step sizes to enhance algorithmic efficiency.

\section{Convergence and Stability analysis}\label{sec:converge}
Based on the results in \cite{Forsythe1968asymptoticdirectionsthes,Yuan2006newstepsizesteepest}, the behavior of the gradient descent method for higher-dimensional problems is essential the same as for two-dimensional problems. Therefore, we consider minimizing the convex quadratic problem \eqref{equ:qua} with  
\begin{equation*}
A=\begin{bmatrix}
\lambda & 0\\
0 & 1
\end{bmatrix},
\end{equation*}
where $\lambda>1$. We denote $\mathbf{g}_{k-1}=(\mathbf{g}_{k-1}^{(1)}, \ \mathbf{g}_{k-1}^{(2)})^{\T}$ with $\mathbf{g}_{k-1}^{(i)}\neq 0$ for $i=1$, $2$ and $k\ge 2$, and define
\begin{equation}\label{equ:epk}
\epsilon_{k}=\frac{(\mathbf{g}_{k-1}^{(1)})^{2}}{(\mathbf{g}_{k-1}^{(2)})^{2}}>0.
\end{equation}
This yields 
\begin{align}\label{equ:normgk}
\begin{split}
\|\g_{k+1}\|_{2}^{2}&=(\g_{k+1}^{(2)})^2(\epsilon_{k+2}+1).
\end{split}	
\end{align}
Let 
\begin{equation*}
\theta_{k}=\angle(\sss_{k-1}, \yy_{k-1}).
\end{equation*}
If $\sss_{k-1}^{\T}\yy_{k-1}>0$, then  
\begin{equation}\label{equ:theta}
\theta_{k}\in[0,\frac{\pi}{2}).
\end{equation}
Define  
\begin{equation*}
p_{k}=1+\sin\theta_{k}.
\end{equation*}
From \eqref{equ:theta}, it follows that 
\begin{equation}
p_{k}\in[1,2).
\end{equation}
Using   $\sss_{k-1}=\x_{k}-\x_{k-1}=-(\alpha_{k-1})^{-1}\g_{k-1}$ and  $\yy_{k-1}=\g_{k}-\g_{k-1}=-(\alpha_{k-1})^{-1}A\g_{k-1}$, we obtain  
\begin{equation}\label{LR}
\alpha_{k}^{L}=\frac{1}{p_{k}}\frac{\sss_{k-1}^{\T}\yy_{k-1}}{\sss_{k-1}^{\T}\sss_{k-1}}=\frac{1}{p_{k}}\frac{\g_{k-1}^{\T}A\g_{k-1}}{\g_{k-1}^{\T}\g_{k-1}},\quad\alpha_{k}^{R}=p_{k}\frac{\yy_{k-1}^{\T}\yy_{k-1}}{\sss_{k-1}^{\T}\yy_{k-1}}=p_{k}\frac{\g_{k-1}^{\T}A^2\g_{k-1}}{\g_{k-1}^{\T}A\g_{k-1}}.
\end{equation}
\begin{property}\label{property1}
	For problem \eqref{pro:qua}, 
	\begin{equation}
	\alpha_{k}^{L}\in(\frac{1}{2}\lambda_{n}, \lambda_{1}],\quad\alpha_{k}^{R}\in[\lambda_{n}, 2\lambda_{1}).
	\end{equation}	
\end{property}
\begin{proof}
	The result follows immediately from  \eqref{A}, \eqref{LR}, and the
	Rayleigh quotient inequality \cite{Golub2013MatrixComputations}.
\end{proof}
\begin{remark}
	Property \ref{property1} is important as it ensures convergence of the LEFT and RIGHT methods.
\end{remark}
From \eqref{equ:epk} and \eqref{LR},  
\begin{equation*}\label{equ:nRBB}
\alpha_{k}^{L}=\frac{1}{p_{k}}\frac{\lambda\epsilon_{k}+1}{\epsilon_{k}+1},\quad\alpha_{k}^{R}=p_{k}\frac{\lambda^2\epsilon_{k}+1}{\lambda\epsilon_{k}+1}.
\end{equation*}
To facilitate analysis, we consider the fixed-parameter variants
\begin{equation}\label{equ:nnRBB}
\alpha_{k}^{L}=\frac{1}{p}\frac{\lambda\epsilon_{k}+1}{\epsilon_{k}+1},\quad\alpha_{k}^{R}=p\frac{\lambda^2\epsilon_{k}+1}{\lambda\epsilon_{k}+1},
\end{equation}
where $p\in[1,2)$ is constant. The following subsections analyze the convergence behavior of the LEFT and RIGHT methods respectively. When $p=1$, these reduce to $\alpha_{k}^{BB1}$ and $\alpha_{k}^{BB2}$. Thus, we restrict our analysis to $p\in(1,2)$.

\subsection{Analysis of the LEFT method}
In the LEFT method, 
\begin{equation*}
\g_{k+1}=\big(I-(\alpha_{k}^{L})^{-1}A\big)\g_{k},
\end{equation*}
where $I$ is the identity matrix. This yields  
\begin{align}
\g_{k+1}^{(1)}&=\big(1-(\alpha_{k}^{L})^{-1}\lambda\big)\g_{k}^{(1)}=\Big(1-\frac{\lambda(\epsilon_{k}+1)p}{\lambda\epsilon_{k}+1
}\Big)\g_{k}^{(1)},\label{equ:gkk11}\\
\g_{k+1}^{(2)}&=\big(1-(\alpha_{k}^{L})^{-1}\big)\g_{k}^{(2)}=\Big(1-\frac{(\epsilon_{k}+1)p}{\lambda\epsilon_{k}+1
}\Big)\g_{k}^{(2)}.\label{equ:gkk22}	
\end{align}
Defining 
\begin{equation}\label{equ:xik}
\xi_{k}=1-(\alpha_{k}^{L})^{-1},
\end{equation}
Property \ref{property1} implies 
\begin{equation*}
(\xi_{k})^2\le 1
\end{equation*}
which indicates that $(\g_{k}^{(2)})^2$ decreases monotonically. By \eqref{equ:normgk}, when $(\xi_{k})^2\le1$, the behavior of  $\|\g_{k+1}\|_{2}$ mirrors that of $\epsilon_{k+2}$: that is, $\|\g_{k+1}\|_{2}$ decrease when $\epsilon_{k+2}\le1$ but may increase otherwise. We therefore focus on the evolution of $\epsilon_{k}$. Define 
\begin{equation}\label{equ:uk}
u(\epsilon_{k})=\frac{1-\lambda[p(\epsilon_{k}+1)-\epsilon_{k}]}{\epsilon_{k}(\lambda-p)-(p-1)}.
\end{equation}
From \eqref{equ:gkk11} and \eqref{equ:gkk22},   
\begin{align}\label{equ:ghk11}
\epsilon_{k+2}=F(\epsilon_{k},\epsilon_{k+1}),
\end{align}
where
\begin{equation}\label{F}
F(\epsilon_{k},\epsilon_{k+1})=g(\epsilon_{k})\epsilon_{k+1},
\end{equation} 
with $g(\epsilon_{k})=(u(\epsilon_{k}))^2$. 

We now analyze LEFT method convergence via $\epsilon_{k}$ dynamics.
\begin{lemma}\label{lemma:lem1}
	If $\epsilon_{k+1}>1$, there exists an integer $N\ge1$ such that
	\begin{equation}
	\epsilon_{k+N+1}\le 1.
	\end{equation} 
\end{lemma}
\begin{proof}
	Assume, for the purpose of deriving a contradiction, that 
	\begin{equation}\label{equ:assume}
	\epsilon_{k+N+1}>1
	\end{equation}
	for all $N\ge1$.
	Denote
	\begin{equation}\label{equ:delta}
	\delta=\frac{1}{\epsilon_{k+N+2}}\in(0,1).
	\end{equation}
	Then 
	\begin{align}\label{equ:large1}
	\begin{split}
	0<\epsilon_{k+N+3}-1&=g(\epsilon_{k+N+1})\epsilon_{k+N+2}-1\\
	&=\big(g(\epsilon_{k+N+1})-\delta\big)\epsilon_{k+N+2}.
	\end{split}
	\end{align}
	Equation \eqref{equ:large1} and the assumption \eqref{equ:assume} imply 
	\begin{equation}\label{equ:lizi2}
	g(\epsilon_{k+N+1})>\delta.
	\end{equation}
	Define  
	\begin{equation}\label{equ:ab}
	a=\frac{(\lambda\epsilon_{k+N+1}+1)(\sqrt{\delta}+1)}{(\epsilon_{k+N+1}+1)(\sqrt{\delta}+\lambda)},\quad b=\frac{(\lambda\epsilon_{k+N+1}+1)(1-\sqrt{\delta})}{(\epsilon_{k+N+1}+1)(\lambda-\sqrt{\delta})}.
	\end{equation}
	Given $\lambda>1$, \eqref{equ:assume}, \eqref{equ:delta}, and \eqref{equ:ab}, 
	\begin{equation}\label{equ:pvalue2}
	a\in(\frac{\lambda+1}{2\lambda},\frac{2\lambda}{\lambda+1})\subset(\frac{1}{2},2),\quad b\in(0,1).
	\end{equation}	
	Condition \eqref{equ:lizi2} requires  
	\begin{equation}\label{equ:pvalue}
	p>a\quad\text{or}\quad p>b,
	\end{equation}	
	corresponding to   
	\begin{equation*}
	\lambda>p+\frac{p-1}{\epsilon_{k+N+1}}\quad\text{or}\quad\lambda<p+\frac{p-1}{\epsilon_{k+N+1}},
	\end{equation*}
	respectively. Since $p\in(1,2)$ cannot satisfy $p>a$ consistently when $\lambda>p+\frac{p-1}{\epsilon_{k+N+1}}$, we obtain a contradiction. 	
\end{proof}
\begin{remark}\label{remark:p2}
	Derivation of \eqref{equ:pvalue2}: For  fixed $\lambda>1$, $a$ increase with $\epsilon_{k+N+1}$ and $\delta$, while $b$  increases with $\epsilon_{k+N+1}$ but  decreases with $\delta$. Lemma \ref{lemma:lem1} implies that if $p\ge2$,  $\epsilon_{k}>1$ for all $k>N$, preventing  $\|\g_{k}\|_{2}$ from decreasing and causing divergence.
\end{remark}
\begin{theorem}\label{thm1}
	For problem \eqref{equ:qua}, let $\{\x_{k}\}$ be generated by the LEFT method. Then either $\g_{k}=\xo$ for some finite $k$, or $\{\|\g_{k}\|_{2}\}$ converges to zero $R$-linearly.
\end{theorem}
\begin{proof}
	Assume $\g_{k}\neq\xo$ for all $k$. If  $\epsilon_{k+1}>1$ occurs infinitely often, $\|\g_{k}\|_{2}\le\|\g_{k-1}\|_{2}$ decreases monotonically. Otherwise, Lemma \ref{lemma:lem1} guarantees an integer $N\ge1$ such that  
	\begin{align}
	\begin{split}
	\|\g_{k+N}\|_{2}^{2}&=(\g_{k+N}^{(2)})^2(1+\varepsilon_{k+N+1})\\
	&=\prod_{j=0}^{N-1}(\xi_{k+j})^2(\g_{k}^{(2)})^2(1+\varepsilon_{k+N+1})\\
	&=\prod_{j=0}^{N-1}(\xi_{k+j})^2\frac{1+\varepsilon_{k+N+1}}{1+\varepsilon_{k+1}}\|\g_{k}\|_{2}^{2}\\
	&<\prod_{j=0}^{N-1}(\xi_{k+j})^2\|\g_{k}\|_{2}^{2},
	\end{split}
	\end{align}
	where $\xi_{k+j}\in(0,1)$. Thus $\{\g_{k}\}$ converges $R$-linearly to zero
\end{proof}

We next analyze the stability of LEFT method, reflected in nonmonotonic oscillations of  $\|\g_{k}\|_{2}$. Certain $p$ values yield stability, revealing spectral gradient step size properties.  

Equation \eqref{equ:ghk11} is a second-order autonomous nonlinear difference equation about $\epsilon_{k}$ \cite[p.2]{Elaydi2005IntroductionDifferenceEquations}. 
Let $\epsilon^{*}$ be an \eqf point satisfying $F(\epsilon^{*},\epsilon^{*})=\epsilon^{*}$, implying    $\epsilon_{k}=\epsilon_{k+1}=\epsilon_{k+2}=\epsilon^{*}$ for $k\ge N$, where $N\ge1$ is an integer.Solving yields:	
\begin{align}\label{equ:fixed}
\epsilon^{*}=0\quad\text{or}\quad\epsilon^{*}=\frac{p(\lambda+1)-2}{\lambda(2-p)-p}.
\end{align}
Key observations:
\begin{enumerate}
	\item[(1)] $\epsilon^{*}=0$ is always an equilibrium point.
	\item[(2)] $\epsilon^{*}=\frac{p(\lambda+1)-2}{\lambda(2-p)-p}$ exits only when $\lambda>\frac{p}{2-p}$.
\end{enumerate}
Linearizing \eqref{equ:ghk11} near $\epsilon^{*}$ via $a_{k}=\epsilon_{k}$, $b_{k}=\epsilon_{k+1}$: 
\begin{equation*}
\begin{cases}
a_{k+1}=b_{k},\\
b_{k+1}=g(a_{k})b_{k},
\end{cases}
\end{equation*}
i.e.,
\begin{equation*}
F\begin{bmatrix}
a_{k}\\
b_{k}
\end{bmatrix}=\begin{bmatrix}
b_{k}\\
g(a_{k})b_{k}
\end{bmatrix}.
\end{equation*}  
The Jacobian is: 
\begin{equation}\label{equ:Jacobian}
J_{1}=\begin{bmatrix}
\frac{\partial F_{1}}{\partial a_{k}},&\frac{\partial F_{1}}{\partial b_{k}}\\
\frac{\partial F_{2}}{\partial a_{k}},&\frac{\partial F_{2}}{\partial b_{k}}
\end{bmatrix}
=\begin{bmatrix}
0,&1\\
b_{k}g'(a_{k}),&g(a_{k})
\end{bmatrix}.
\end{equation}

\textbf{Case 1}:
$$\epsilon^{*}=0.$$ 
At $(a_{k},b_{k})^{\T}=(0,0)^{\T}$, 
$g(0)=(\frac{\lambda p-1}{p-1})^2>1$ 
($\lambda>1$, $p\in(1,2)$). Eigenvalues of 
$J_{1}$ are $0$ and $g(0)$, implying instability (\textbf{Theorem1}).

\textbf{Case 2}:
\begin{equation*}\label{equ:ekk}
\epsilon^{*}=\frac{p(\lambda+1)-2}{\lambda(2-p)-p}.
\end{equation*}
Here $u(\epsilon^{*})=-1$, $g(\epsilon^{*})=1$. The Jacobian at $(\epsilon^{*},\epsilon^{*})^{\T}$ is
\begin{equation*}\label{equ:Jacobian2}
J_{1}=\begin{bmatrix}
0,&1\\
\epsilon^{*}g'(\epsilon^{*}),&1
\end{bmatrix},
\end{equation*}
with characteristic equation: 
\begin{equation}\label{equ:cha}
\mu^{2}-\mu+q=0,
\end{equation}
where
\begin{equation}\label{equ:q}
q=-\epsilon^{*}g'(\epsilon^{*})=\frac{2[p(\lambda+1)-2][\lambda(2-p)-p]}{p(\lambda-1)^2}\in(-\infty,2).	
\end{equation}
Stability depends on the discriminant $D=1-4q:$ 
\begin{enumerate}
	\item[(1)] $D>0$, i.e., $q<\frac{1}{4}$. Roots 
	\begin{equation*}
	\mu=\frac{1\pm\sqrt{1-4q}}{2}.
	\end{equation*}
	If $q\in(0,\frac{1}{4})$, then we have
	\begin{equation*}
	\mu_{1}\in(0,0.5),\quad\mu_{2}\in(0.5,1).
	\end{equation*}
	Both solutions satisfy $|\mu|<1$. Thus, the state $(\epsilon^{*},\epsilon^{*})^{\T}$ is stable in this scenario.
	\item[(2)] $D=0$, i.e., $q=\frac{1}{4}$. Then we have
	\begin{equation*}
	\mu=\frac{1}{2}.
	\end{equation*}
	We know that the state $(\epsilon^{*},\epsilon^{*})^{\T}$ is stable. 
	\item[(3)] $D<0$, i.e., $q>\frac{1}{4}$. The solutions to \eqref{equ:cha} are
	\begin{equation}\label{equ:complexmu}
	\mu=\frac{1\pm i\sqrt{4q-1}}{2}.
	\end{equation} 
	From \eqref{equ:complexmu}, we have 
	\begin{equation}
	|\mu|=\sqrt{(\frac{1}{2})^2+(\frac{\sqrt{4q-1}}{2})^2}=\sqrt{q}.
	\end{equation}
	Thus, the state $(\epsilon^{*},\epsilon^{*})^{\T}$ is stable if $q\in(\frac{1}{4},1)$. When $q\in(1,2)$, then the state $(\epsilon^{*},\epsilon^{*})^{\T}$ is unstable. 
\end{enumerate}

\begin{figure}[!h]	
	\centering
	\subfigure{
		\includegraphics[width=0.9\textwidth]{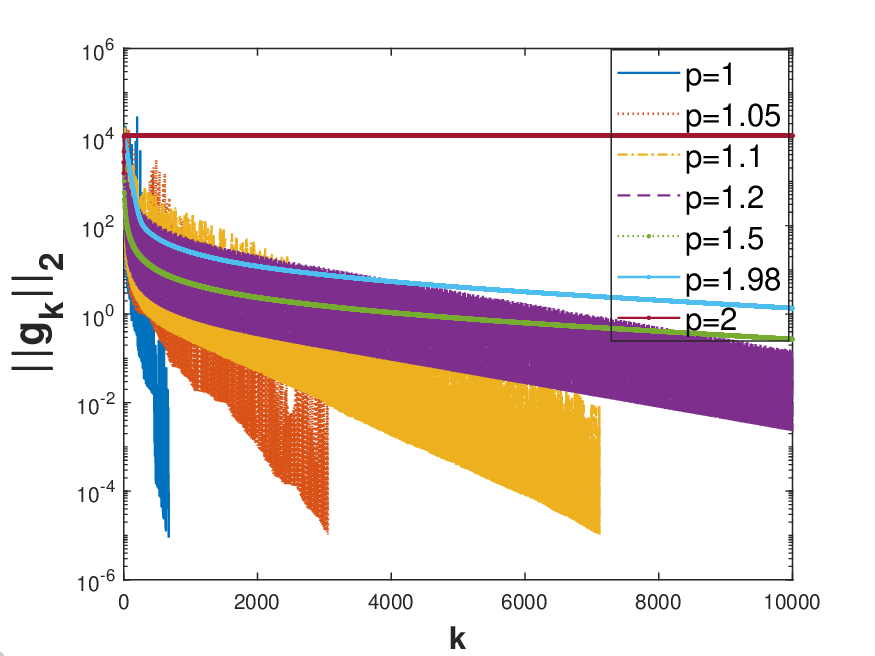}}\\	
	\caption{\textit{The performance of the LEFT method with different $p$ in problem \eqref{pro:qua}.}}	
	\label{fig:BB1V}
\end{figure}

Based on these analyses, we conclude as follows.
\begin{theorem}
	For \eqf point $\epsilon^{*}$ of \eqref{equ:ghk11}:
	\begin{enumerate}	
		\item If $\lambda>\frac{p}{2-p}$ and $q\in(0,1)$, then there exists a stable positive $\epsilon^{*}$. 
		\item If $\lambda>\frac{p}{2-p}$ and $q\in(-\infty,0)\cup(1,2)$, there exists a positive $\epsilon^{*}$, but it is unstable, and the sequence $\{\epsilon_{k}\}$ may oscillate wildly.
		\item If $\lambda\le\frac{p}{2-p}$,  $\epsilon^{*}=0$ is unstable, and the sequence $\{\epsilon_{k}\}$ may oscillate. 
	\end{enumerate}
\end{theorem}
From \eqref{equ:q}, we know that $q\rightarrow2(2-p)$ as $\lambda\rightarrow\infty$ (i.e., $q\in(0,2)$). Based on this, if $p\in[1.5,2)$, then $q\in(0,1)$, that is, the \eqf point $\epsilon^{*}$ is stable. Conversely, if $p\in[1,1.5)$, then the \eqf point $\epsilon^{*}$ may be unstable. We test the effect of the p-value on the performance of the LEFT method in problem \eqref{pro:qua}, using parameter values $p=1$, $1.05$, $1.1$, $1.2$, $1.5$, $1.98$, and $2$. Figure \ref{fig:BB1V} shows the performance of LEFT method with these $p$ values. Note that the LEFT method with $p=1$ corresponds the BB1 method.

Figure \ref{fig:BB1V} validates our analysis: when $p\in[1,1.5)$, $\|\g_{k}\|_{2}$ exhibits violent oscillation; when $p\in[1.5,2)$, $\|\g_{k}\|_{2}$ demonstrates a stable convergence pattern. Notably, when $p=2$, $\|\g_{k}\|_{2}$ fails to convergence. This indicates that $\alpha_{k}^{L}$ in \eqref{equ:BB1left} has reached the boundary of convergence--that is, selecting smaller scalar to generate the step size would cause the algorithm to diverge. This numerical phenomenon verifies the conclusion of Remark \ref{remark:p2}.
\begin{remark}
	From \eqref{equ:fixed}, we know that the positive fixed point $\epsilon^{*}\rightarrow\frac{p}{2-p}$ as $\lambda\rightarrow\infty$. Based on this, we derive that $\epsilon^{*}\rightarrow1$ as $p\rightarrow1$, and $\epsilon^{*}\rightarrow\infty$ as $p\rightarrow2$, that is, the asymptotic convergence rate of the LEFT method decreases with the increase of $p$, while its stability increases with the increase of $p$. From this, we get a derivative result: the asymptotic convergence rate of the original BB1 method is approximately $\xi_{k}$ \eqref{equ:xik} due to $\epsilon^{*}\rightarrow1$.
\end{remark}

\subsection{Analysis of the RIGHT method}
In the RIGHT method, 
\begin{equation*}
\g_{k+1}=\big(I-(\alpha_{k}^{R})^{-1}A\big)\g_{k}.
\end{equation*}
Then we have 
\begin{align}
\g_{k+1}^{(1)}&=\big(1-(\alpha_{k}^{R})^{-1}\lambda\big)\g_{k}^{(1)}=\Big(1-\frac{\lambda(\lambda\epsilon_{k}+1)}{p(\lambda^2\epsilon_{k}+1)}\Big)\g_{k}^{(1)},\label{equ:gkk211}\\
\g_{k+1}^{(2)}&=\big(1-(\alpha_{k}^{R})^{-1}\big)\g_{k}^{(2)}=\Big(1-\frac{(\lambda\epsilon_{k}+1)}{p(\lambda^2\epsilon_{k}+1)}\Big)\g_{k}^{(2)}.\label{equ:gkk222}	
\end{align}
Define
\begin{equation}\label{equ:eta}
\eta_{k}=1-\frac{(\lambda\epsilon_{k}+1)}{p(\lambda^2\epsilon_{k}+1)}.
\end{equation}
This implies  
\begin{equation*}
\eta_{k}\le 1.
\end{equation*}
By \eqref{equ:normgk}, the behavior of  $\|\g_{k+1}\|_{2}$ is determined by  $\epsilon_{k+2}$. We therefore analyze the evolution of $\epsilon_{k}$. Define  
\begin{equation}\label{equ:vk}
v(\epsilon_{k})=\frac{\lambda^2\epsilon_{k}(p-1)+p-\lambda}{\lambda\epsilon_{k}(p\lambda-1)+p-1}.
\end{equation}
From \eqref{equ:gkk211} and \eqref{equ:gkk222},    
\begin{align}\label{equ:ghk211}
\epsilon_{k+2}=G(\epsilon_{k},\epsilon_{k+1}).
\end{align}
with 
\begin{equation*}
G(\epsilon_{k},\epsilon_{k+1})=l(\epsilon_{k})\epsilon_{k+1},
\end{equation*}
where $l(\epsilon_{k})=(v(\epsilon_{k}))^2$.

The convergence analysis for the RIGHT method follows a similar approach to the LEFT method.
\begin{lemma}\label{lemma:lem2}
	If $\epsilon_{k+1}>1$, then there exists an integer $N\ge1$ such that
	\begin{equation}
	\epsilon_{k+N+1}\le 1.
	\end{equation} 
\end{lemma}
\begin{proof}
	Assume by contradiction that 
	\begin{equation}\label{equ:assume2}
	\epsilon_{k+N+1}>1
	\end{equation}
	for all $N\ge1$. Let 
	\begin{equation}\label{equ:delta2}
	\delta=\frac{1}{\epsilon_{k+N+2}}\in(0,1).
	\end{equation}
	Then 
	\begin{align}\label{equ:large21}
	\begin{split}
	0<\epsilon_{k+N+3}-1&=l(\epsilon_{k+N+1})\epsilon_{k+N+2}-1\\
	&=\big(l(\epsilon_{k+N+1})-\delta\big)\epsilon_{k+N+2}.
	\end{split}
	\end{align}
	It follows from \eqref{equ:assume2} and \eqref{equ:large21} that 
	\begin{equation}\label{equ:lizi3}
	l(\epsilon_{k+N+1})>\delta.
	\end{equation}
	Define  
	\begin{equation}\label{equ:ab2}
	a=\frac{(\lambda\epsilon_{k+N+1}+1)(\lambda-\sqrt{\delta})}{(\lambda^2\epsilon_{k+N+1}+1)(1-\sqrt{\delta})},\quad b=\frac{(\lambda\epsilon_{k+N+1}+1)(\lambda+\sqrt{\delta})}{(\lambda^2\epsilon_{k+N+1}+1)(1+\sqrt{\delta})}.
	\end{equation}
	Given $\lambda>1$, \eqref{equ:assume2}, and  \eqref{equ:delta2}, we have
	\begin{equation}\label{equ:pvalue4}
	a\in(1,\infty)\quad\text{and}\quad b\in(\frac{1}{2},\frac{4+3\sqrt{2}}{4+2\sqrt{2}}).
	\end{equation}	
	Condition \eqref{equ:lizi3} requires 
	\begin{equation}\label{equ:pvalue3}
	p>a\quad\text{or}\quad p<b,
	\end{equation}	
	corresponding to 
	\begin{equation*}
	(\lambda^2\epsilon_{k+N+1}+1)p-\lambda^2\epsilon_{k+N+1}-\lambda>0\quad\text{or}\quad(\lambda^2\epsilon_{k+N+1}+1)p-\lambda^2\epsilon_{k+N+1}-\lambda<0,
	\end{equation*}
	respectively. Since $p\in(1,2)$ cannot satisfy $p>a$ consistently when $(\lambda^2\epsilon_{k+N+1}+1)p-\lambda^2\epsilon_{k+N+1}-\lambda>0$, we obtain a contradiction.
\end{proof}
\begin{remark}\label{remark:p3}
	The derivation of \eqref{equ:pvalue4} comes from the following fact: For a fixed $\lambda>1$, $b$ in \eqref{equ:ab2} is monotonically decreasing with respect to both $\epsilon_{k+N+1}$ and $\delta$, while $a$ is monotonically increasing with respect to $\delta$ and monotonically decreasing with respect to $\epsilon_{k+N+1}$. Unlike the LEFT method, Lemma \ref{lemma:lem2} shows that no matter what value $p>1$ takes, $\epsilon_{k}>1$ cannot be always true. In other words, for any $p>1$, the RIGHT method always converges.
\end{remark}
\begin{theorem}\label{thm2}
	For problem \eqref{equ:qua}, let $\{\x_{k}\}$ be generated by the RIGHT method. Then either $\g_{k}=\xo$ for some finite $k$, or $\{\|\g_{k}\|_{2}\}$ converges to zero $R$-linearly.
\end{theorem}
\begin{proof}
	We only need to consider the case that $\g_{k}\neq\xo$ for all $k$. We still assume that $\epsilon_{k+1}>1$ for some $k>1$. From \eqref{equ:normgk}, \eqref{equ:gkk222}, and Lemma \ref{lemma:lem2} we know that there exists an integer $N\ge1$ such that  
	\begin{align}
	\begin{split}
	\|\g_{k+N}\|_{2}^{2}&=(\g_{k+N}^{(2)})^2(1+\epsilon_{k+N+1})\\
	&=\prod_{j=0}^{N-1}(\eta_{k+j})^2(\g_{k}^{(2)})^2(1+\epsilon_{k+N+1})\\
	&=\prod_{j=0}^{N-1}(\eta_{k+j})^2\frac{1+\epsilon_{k+N+1}}{1+\epsilon_{k+1}}\|\g_{k}\|_{2}^{2}\\
	&<\prod_{j=0}^{N-1}(\eta_{k+j})^2\|\g_{k}\|_{2}^{2},
	\end{split}
	\end{align}
	where each $\eta_{k+j}\in(0,1)$. Therefore  $\{\g_{k}\}$ converges to zero $R$-linearly.	
\end{proof}

We now analyze the stability of the RIGHT method. Similar to \eqref{equ:ghk11}, equation  \eqref{equ:ghk211} is also a second-order autonomous nonlinear system of $\epsilon_{k}$. We examine its the \eqf points and their stability. Let $\epsilon^{*}$ satisfy  $\epsilon^{*}=G(\epsilon^{*},\epsilon^{*})$, and solve for $\epsilon^{*}$, implying  $\epsilon_{k}=\epsilon_{k+1}=\epsilon_{k+2}=\epsilon^{*}$ for all $k\ge N$, where $N\ge1$ is an integer. Solving yields: 	
\begin{align}\label{equ:fixed2}
\epsilon^{*}=0,\quad\quad\epsilon^{*}=\frac{\lambda+1-2p}{\lambda[\lambda(2p-1)-1]},\quad\text{or}\quad\epsilon^{*}=-\frac{1}{\lambda}.
\end{align}
Since $\epsilon^{*}=-\frac{1}{\lambda}$ contradicts $\epsilon_{k}\ge0$, we discard it. We consider: 
\begin{enumerate}
	\item $\epsilon^{*}=0$ is an \eqf point for any $\lambda>1$.
	\item $\epsilon=\frac{\lambda+1-2p}{\lambda[\lambda(2p-1)-1]}$ is an \eqf point only when $\lambda\ge2p-1$. 
\end{enumerate}

Linearizing \eqref{equ:ghk211} near $\epsilon^{*}$ via $c_{k}=\epsilon_{k}$, $d_{k}=\epsilon_{k+1}$:
\begin{equation}
\begin{cases}
c_{k+1}=d_{k},\\
d_{k+1}=l(c_{k})d_{k}.
\end{cases}
\end{equation}
i.e.,
\begin{equation*}
G\begin{bmatrix}
c_{k}\\
d_{k}
\end{bmatrix}=\begin{bmatrix}
d_{k}\\
l(c_{k})d_{k}
\end{bmatrix}.
\end{equation*}  
Then its Jacobian matrix at $(c_{k},d_{k})^{\T}$ is 
\begin{equation}\label{equ:Jacobia2}
J_{2}=\begin{bmatrix}
\frac{\partial G_{1}}{\partial c_{k}},&\frac{\partial G_{1}}{\partial d_{k}}\\
\frac{\partial G_{2}}{\partial c_{k}},&\frac{\partial G_{2}}{\partial d_{k}}
\end{bmatrix}
=\begin{bmatrix}
0,&1\\
d_{k}l'(c_{k}),&l(c_{k})
\end{bmatrix}.
\end{equation}
We first consider the \eqf point 
$$\epsilon^{*}=0.$$ 
In this case, the state vector $(\epsilon^{*},\epsilon^{*})^{\T}=(0,0)^{\T}$, and the corresponding $l(\epsilon^{*})=(\frac{p-\lambda}{p-1})^2$. The eigenvalues of the matrix $J_{2}$ in \eqref{equ:Jacobia2} are $\mu_{1}=0$, and $\mu_{2}=l(\epsilon^{*})$.
From Theorem \ref{thm:stab2}, we know that if
\begin{equation*}\label{equ:s1}
(\frac{p-\lambda}{p-1})^2<1,
\end{equation*} 
then the \eqf point $\epsilon^{*}=0$ is stable, i.e., 
\begin{equation*}
\lambda\in\big(1,2p-1\big).
\end{equation*}
Therefore, when $\lambda\ge2p-1$, the \eqf point $\epsilon^{*}=0$ is unstable.

We next consider the \eqf point 
\begin{equation}\label{equ:fix}
\epsilon^{*}=\frac{\lambda+1-2p}{\lambda[\lambda(2p-1)-1]}.
\end{equation}
In this case, $v(\epsilon^{*})=-1$, and $l(\epsilon^{*})=1$. Then the Jacobian matrix
\begin{equation*}
J_{2}=\begin{bmatrix}
0,&1\\
-2\epsilon^{*}v'(\epsilon^{*}),&1
\end{bmatrix},
\end{equation*}
and its corresponding characteristic equation is 
\begin{equation}\label{equ:char}
\mu^{2}-\mu+2\epsilon^{*}v'(\epsilon^{*})=0.
\end{equation}
The solutions to \eqref{equ:char} are
\begin{equation}\label{equ:mu12}
\mu=\frac{1\pm\sqrt{1-8\epsilon^{*} v'(\epsilon^{*})}}{2},
\end{equation}
where 
\begin{equation*}\label{equ:vk'}
v'(\epsilon^{*})=\frac{p\lambda(\lambda-1)^2}{(c\epsilon^{*}+d)^2}>0,
\end{equation*}
where $c=\lambda(p\lambda-1)>0$, $d=p-1>0$. Let 
\begin{equation}\label{equ:gamma}
\gamma=\epsilon^{*} v'(\epsilon^{*})=\frac{(\lambda+1-2p)(\lambda(2p-1)-1)}{p(\lambda-1)^2}>0.
\end{equation}
From the relationship between the roots and coefficients of quadratic equations, we know that if the equation \eqref{equ:char} has two solutions $\mu_{1}$ and $\mu_{2}$, then $\mu_{1}>0$, $\mu_{2}>0$. Let $\mu_{1}<\mu_{2}$. Given that the discriminant $D=1-8\gamma$ in \eqref{equ:char}, we consider the following two cases.
\begin{enumerate}
	\item $D>0$, i.e., $\gamma\in(0,\frac{1}{8})$. Because $\mu_{2}=\frac{1+\sqrt{D}}{2}<1$, and $\mu_{1}\le\mu_{2}$, the \eqf point $\epsilon^{*}$ is stable.
	\item $D<0$, i.e., $\gamma>\frac{1}{8}$. In this case, $\mu=\frac{1\pm i\sqrt{|D|}}{2}$. Based on this, from \textbf{Theorem}, if $\gamma\in(\frac{1}{8},\frac{1}{2})$, then the \eqf point $\epsilon^{*}$ is stable.
\end{enumerate}
Thus, if $\gamma\in(0,\frac{1}{2})$, then the state $(\epsilon^{*},\epsilon^{*})^{\T}$ is stable. Nevertheless, from \eqref{equ:gamma}, when $\lambda>2+\sqrt{3}$, $\gamma<0.5$ has no solution in $p\in(1,2)$. This shows that when $\lambda>2+\sqrt{3}$, the \eqf point $\epsilon^{*}$ is unstable for any $p\in(1,2)$. 

\begin{remark}
	Although the \eqf point $\epsilon^{*}=0$  and $\epsilon^{*}=\frac{\lambda+1-2p}{\lambda[\lambda(2p-1)-1]}$ in the RIGHT method are unstable for large $\lambda$, we know that the \eqf point $\epsilon^{*}=\frac{\lambda+1-2p}{\lambda[\lambda(2p-1)-1]}\rightarrow0$ as $\lambda\rightarrow\infty$. Therefore, for $p\in(1,2)$, the global convergence of the RIGHT method is favorable, particular for some ill-conditioned problems. On the other hand, the convergence rate $\eta_{k}$ \eqref{equ:eta} of the RIGHT method approximates $1-\frac{1}{p\lambda}$ as $\lambda\rightarrow\infty$. This indicates that as $p$ increases from $1$ to $2$, $\eta_{k}$ rises, meaning that within the RIGHT method, the BB2 method corresponding to $p=1$ exhibits the best convergence.   
\end{remark}

\begin{figure}[h!]
	\centering
	\subfigure{\includegraphics[width=0.9\linewidth]{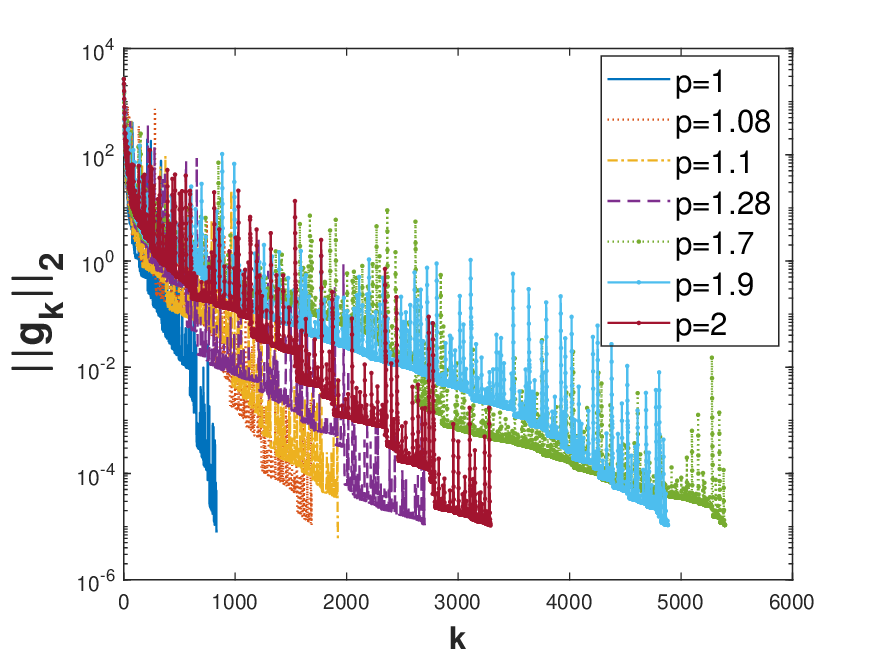}}
	\caption{\textit{The performance of the RIGHT method with different $p$ in  problem \eqref{pro:qua}.}}
	\label{fig:BB2V}
\end{figure}

We test the RIGHT method on problem \eqref{pro:qua} with  $p=1$, $1.08$, $1.1$, $1.28$, $1.7$, $1.9$, and $2$. Figure \ref{fig:BB2V} shows the performance of the RIGHT method with these $p$ values. The results in Figure \ref{fig:BB2V} verify the this conclusion, that is, for large $\lambda$, the RIGHT method corresponding to $p\in(1,2)$ is unstable.

\section{Two modified step sizes}\label{sec:twomodify}
In Section \ref{sec:converge}, we analyze the convergence behavior of the LEFT and RIGHT methods from a dynamics perspective. We find that although both methods converge, their asymptotic convergence rates are slower than those of the original two BB methods, and this characteristic is confirmed by numerical results. In this section, we analyze the underlying reasons for this phenomenon and, based on these findings, propose two modified step sizes. 

Based on the analysis in Section \ref{sec:converge}, we know that the convergence rates of the LEFT and RIGHT methods are primarily determined by $\xi_{k}$ \eqref{equ:xik} and $\eta_{k}$ \eqref{equ:eta}, respectively. From Property \ref{property1}, we have
\begin{equation}\label{equ:xxik}
\xi_{k}\in(-1,0]\cup[0,1-\frac{1}{\lambda}],
\end{equation}
and
\begin{equation}\label{equ:eetak}
\eta_{k}\in[0,1-\frac{1}{\lambda}]\cup[1-\frac{1}{\lambda},1-\frac{1}{2\lambda}).
\end{equation}
Observe \eqref{equ:xxik}, for $(\xi_{k})^{2}$, values in $(-1,0]$ and $[0,1-\frac{1}{\lambda}]$ have a similar effect, causing an increase in the algorithm's iteration number and resulting in a decrease in convergence speed. For $(\eta_{k})^2$, the presence of $[1-\frac{1}{\lambda},1-\frac{1}{2\lambda})$ brings it closer to $1$, which slows down its convergence rate. The fundamental reason for these phenomena is that the ranges of $\alpha_{k}^{L}$ and $\alpha_{k}^{R}$ exceed the range of the eigenvalues of the Hessian matrix $A$, leading to a degradation in the algorithm's convergence rate. 

The philosophy behind the LEFT and RIGHT methods lies in the fact that, compared to the original BB1 and BB2 methods, the scalars $\alpha_{k}^{L}$ and $\alpha_{k}^{R}$ they generate can rapidly approximate the smallest and largest eigenvalues of Hessian $A$, respectively. This ability to quickly scan the spectrum of $A$ is crucial to the performance of spectral gradient descent methods, with relevant theories and algorithms detailed in references such as \cite{Frassoldati2008Newadaptivestepsize,Ferrandi2023harmonicframeworkstepsize,DiSerafino2018steplengthselectiongradient,An2025RegularizedBarzilaiBorwein}. Nevertheless, as elaborated in the preceding paragraph, $\alpha_{k}^{L}$ and $\alpha_{k}^{R}$ may be smaller than the smallest eigenvalue and larger than the largest eigenvalue of Hessian $A$ respectively, which in turn affects the algorithm's convergence rate.

A feasible remedial measure is impose constraints on $\alpha_{k}^{L}$ and $\alpha_{k}^{R}$, ensuring they remain within the spectral radius of $A$ as follows
\begin{equation}
\alpha_{k}^{ML}=\max\{\alpha_{k-1}^{BB1},\alpha_{k}^{L}\}=\begin{cases}
\alpha_{k-1}^{BB1},&\quad\text{if}\quad\alpha_{k-1}^{BB1}\ge\alpha_{k}^{L},\\
\alpha_{k}^{L},&\quad\text{otherwise},
\end{cases}
\end{equation}
\begin{equation}
\alpha_{k}^{MR}=\min\{\alpha_{k-1}^{BB2},\alpha_{k}^{R}\}=\begin{cases}
\alpha_{k-1}^{BB2},&\quad\text{if}\quad\alpha_{k-1}^{BB2}\le\alpha_{k}^{R},\\
\alpha_{k}^{R},&\quad\text{otheriwise}.
\end{cases}
\end{equation}
We term the gradient descent methods with step sizes $\frac{1}{\alpha_{k}^{ML}}$ and $\frac{1}{\alpha_{k}^{MR}}$ the Modified LEFT (ML) and Modified RIGHT (MR)methods, respectively.

\begin{figure}[h!]
	\centering
	\subfigure{\includegraphics[width=0.9\linewidth]{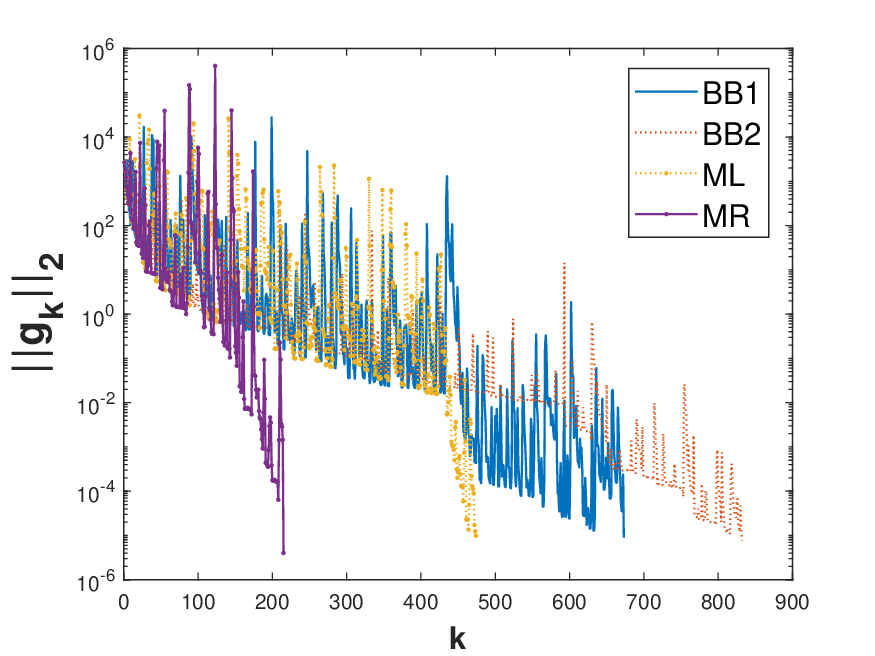}}
	\caption{\textit{The performance of the BB1, BB2, LEFT, RIGHT, LM, and RM methods in the problem \eqref{pro:qua}, with $n=10$.}}
	\label{fig:LRMGnorm}
\end{figure}

We evaluate the ML and MR methods again on problem \eqref{pro:qua}, using the same initial conditions and stopping criteria as in Section \ref{sec:twostep}. Figure \ref{fig:LRMGnorm} shows the performance of BB1, BB2, ML, and MR methods, measured by the Euclidean $2$-norm of the gradient. 

As shown in Figure \ref{fig:LRMGnorm}, the MR method exhibits near-superlinear convergence in this small-scale problem, while the ML method also outperforms both the BB1 and BB2 methods.

\section{Numerical experiments}\label{sec:numerical}
In this section, we present preliminary numerical experiments\footnote{All experiments were implemented in \text{MATLAB R2024a}. All the runs were carried out on a PC with an 12th Gen Intel(R) Core(TM) i7-12700H 2.30 GHz and 32 GB of RAM.} to evaluate the proposed methods.  Performance is assessed using performance profiles \cite{Dolan2002Benchmarkingoptimizationsoftware}, where the solution cost for each problem is normalized by the lowest cost to obtain the performance ratio $\tau$. The most efficient method solves a given problem with a performance ratio $1$, while all other methods solve the problem with a performance ratio of at least $1$. In order to grasp the full implications of our test data regarding the solvers' probability of successfully handing a problem, we display a $log$ scale of the performance profiles. Since we are also interested in the behavior for $\tau$ close to $1$, we use use a base of $2$ for the scale \cite{Dolan2002Benchmarkingoptimizationsoftware}. Therefore, the value of $\rho_{s}(\log_2(1))$ is the probability that the solver $s$ will win over the rest of the solvers. Unless otherwise specified, the performance profiles mentioned in subsequent experiments are all $\log_2$ scaled. For convenience, we denote $\omega=\log_2(\tau)$ in this paper.

\begin{figure}[!h]
	\centering
	\subfigure{
		\includegraphics[width=0.33\textwidth]{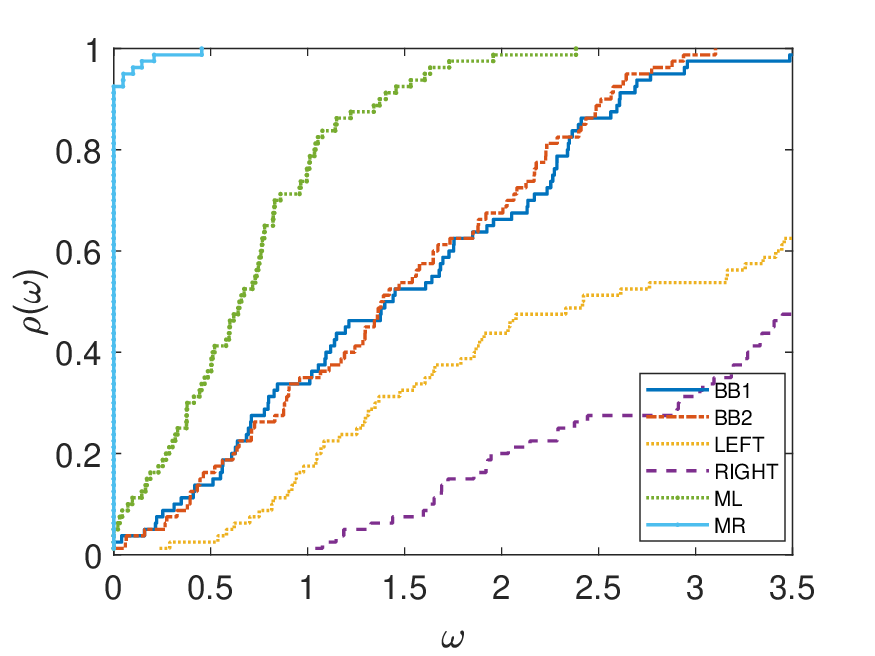}}\hspace{-14pt}
	\subfigure{
		\includegraphics[width=0.33\textwidth]{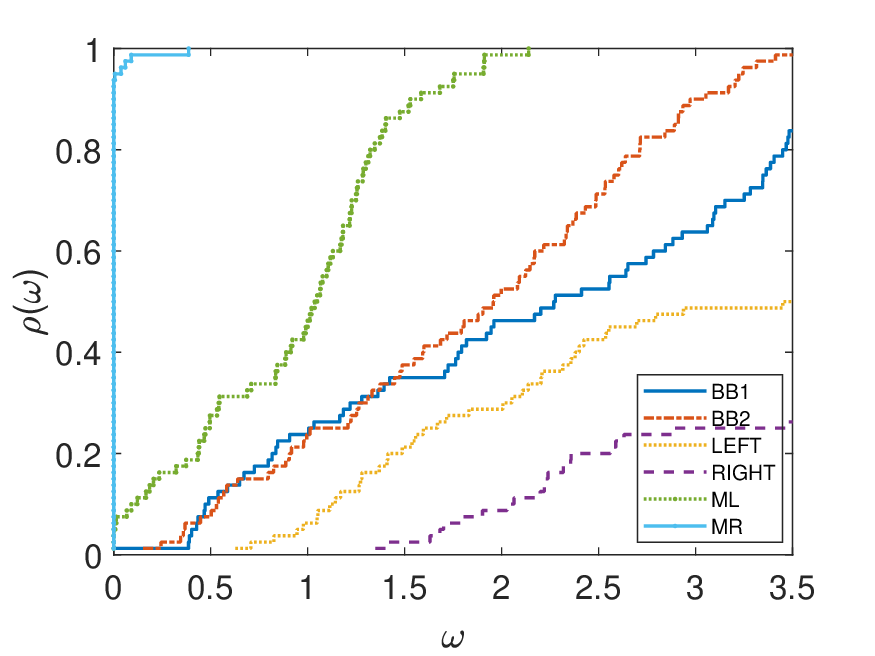}}\hspace{-14pt}
	\subfigure{
		\includegraphics[width=0.33\textwidth]{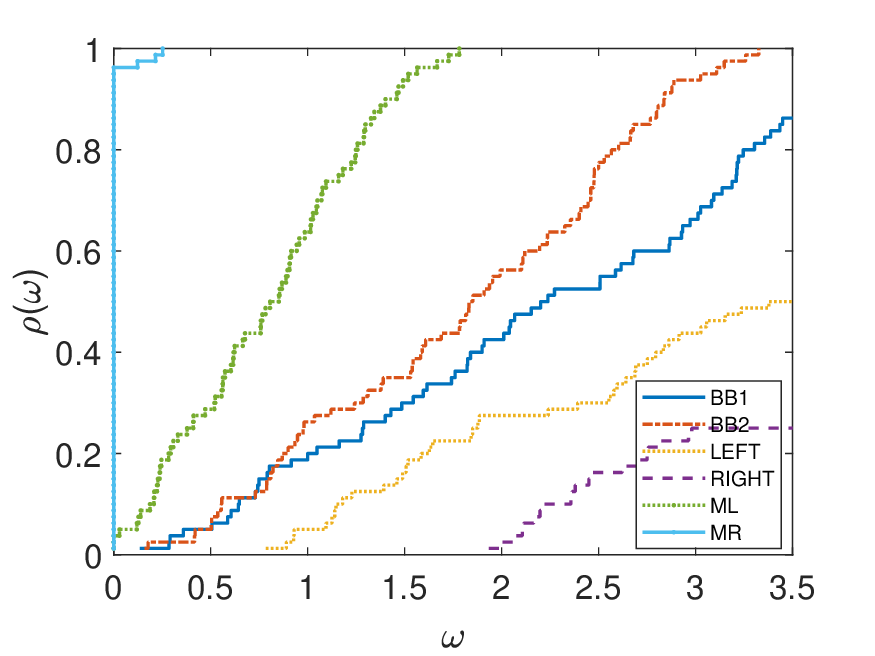}}\\\vspace{-12pt}
	\subfigure{
		\includegraphics[width=0.33\textwidth]{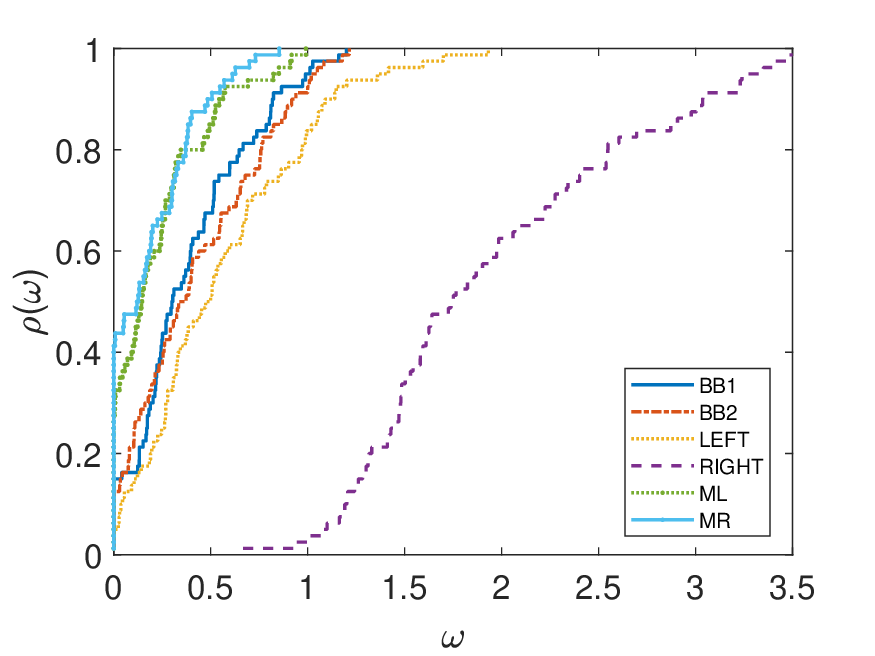}}\hspace{-14pt}
	\subfigure{
		\includegraphics[width=0.33\textwidth]{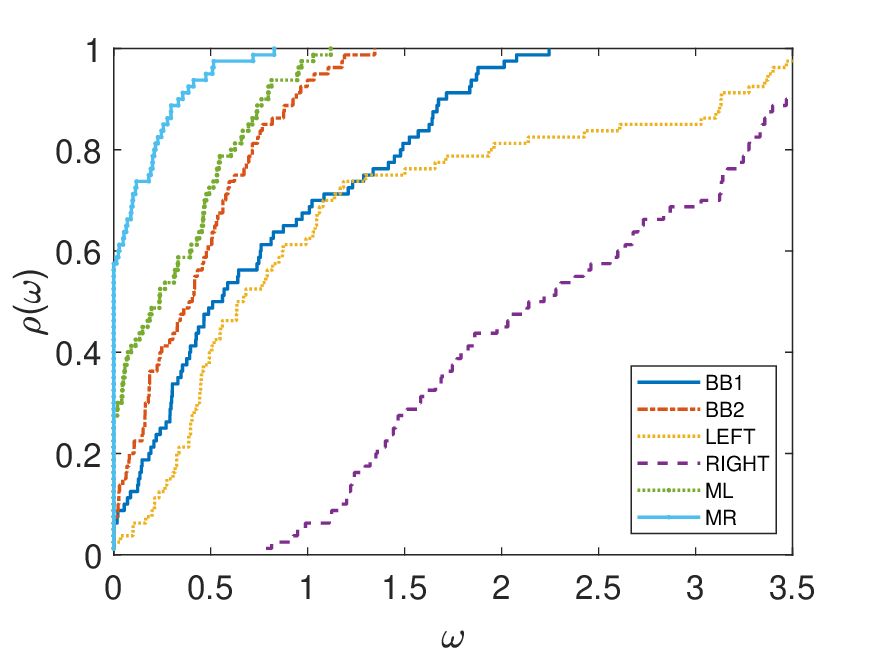}}\hspace{-14pt}
	\subfigure{
		\includegraphics[width=0.33\textwidth]{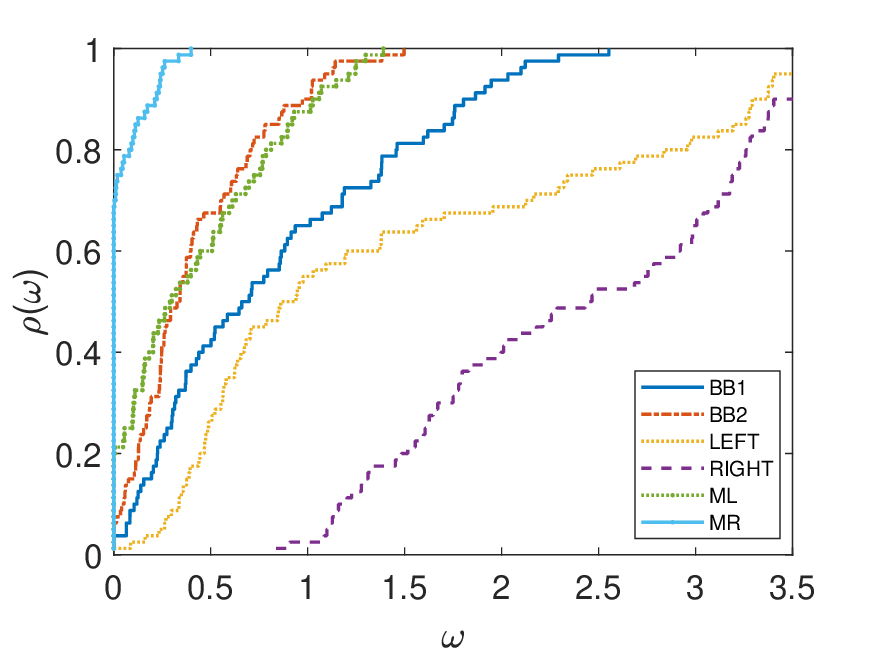}}\\\vspace{-12pt}
	\subfigure{
		\includegraphics[width=0.33\textwidth]{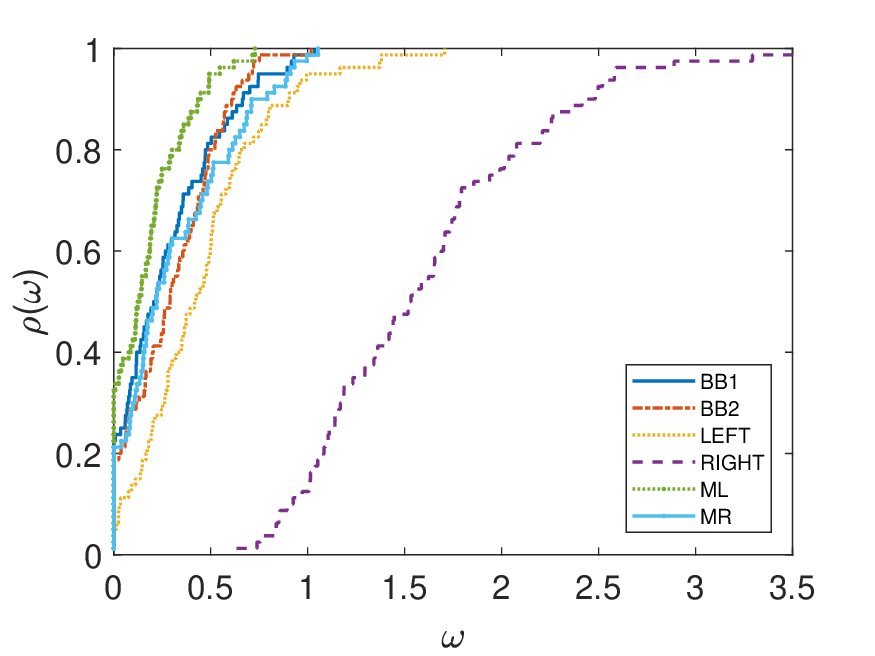}}\hspace{-14pt}
	\subfigure{
		\includegraphics[width=0.33\textwidth]{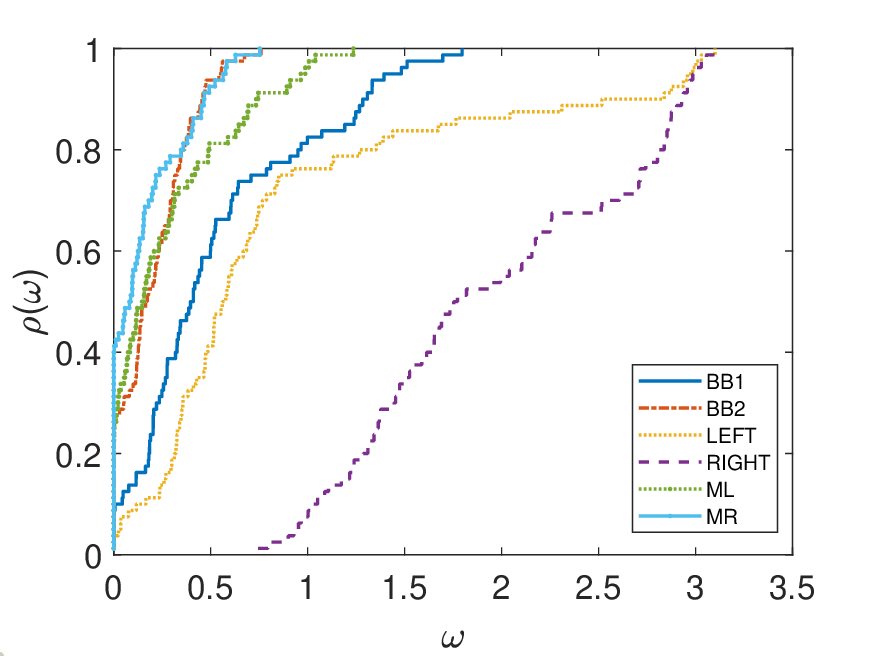}}\hspace{-14pt}
	\subfigure{
		\includegraphics[width=0.33\textwidth]{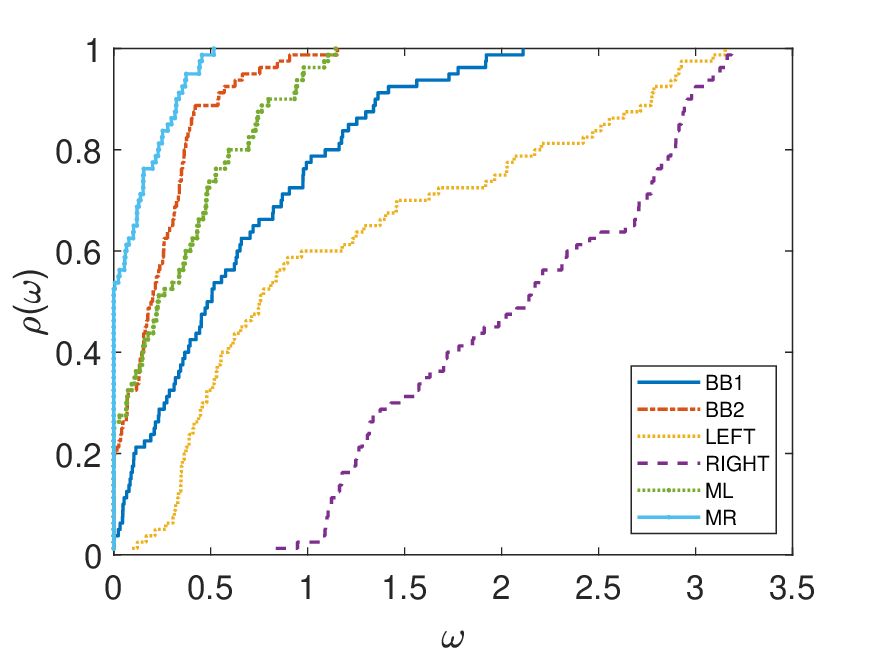}}\\
	\caption{\textit{Performance profiles of the these method on the quadratic problems \eqref{pro:qua}, from the first row to the last row, $n=10,50,100$, from the first column to the last column, $\varepsilon=10^{-6},10^{-9},10^{-12}$, based on iteration number.}}	
	\label{fig:nsq1}
\end{figure}

\begin{figure}[!h]
	\centering	
	\subfigure{
		\includegraphics[width=0.33\textwidth]{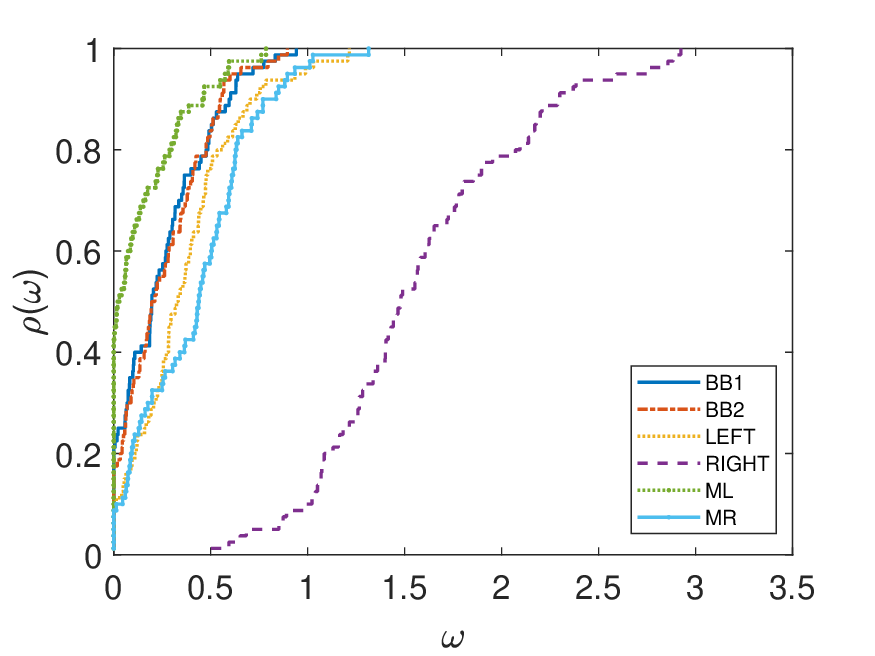}}\hspace{-14pt}
	\subfigure{
		\includegraphics[width=0.33\textwidth]{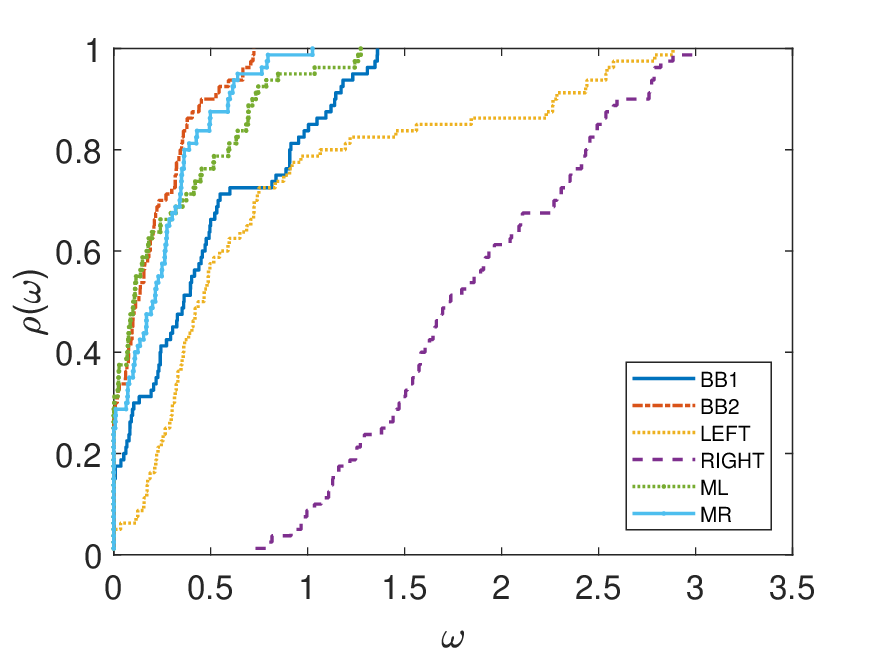}}\hspace{-14pt}
	\subfigure{
		\includegraphics[width=0.33\textwidth]{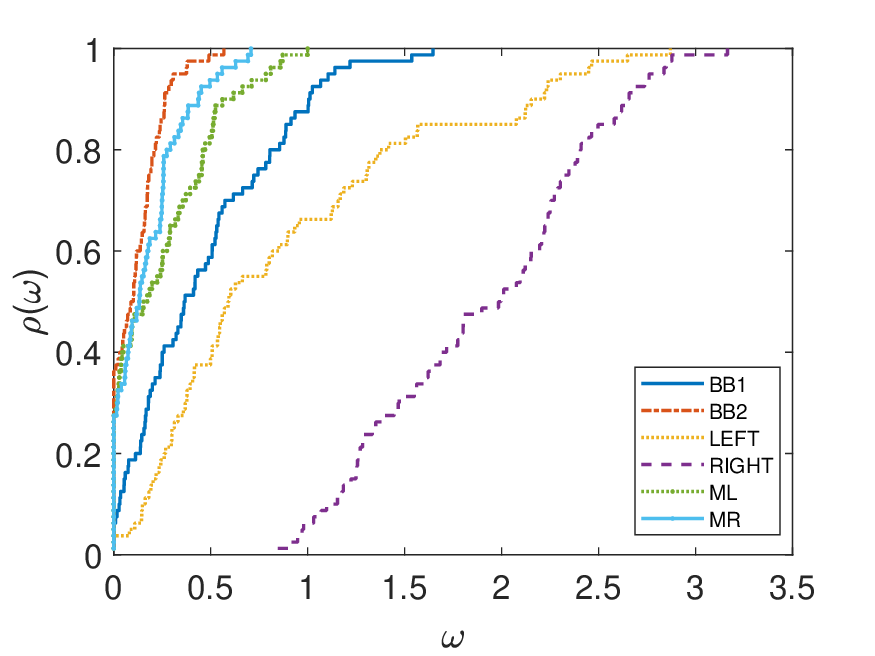}}\\\vspace{-12pt}
	\subfigure{
		\includegraphics[width=0.33\textwidth]{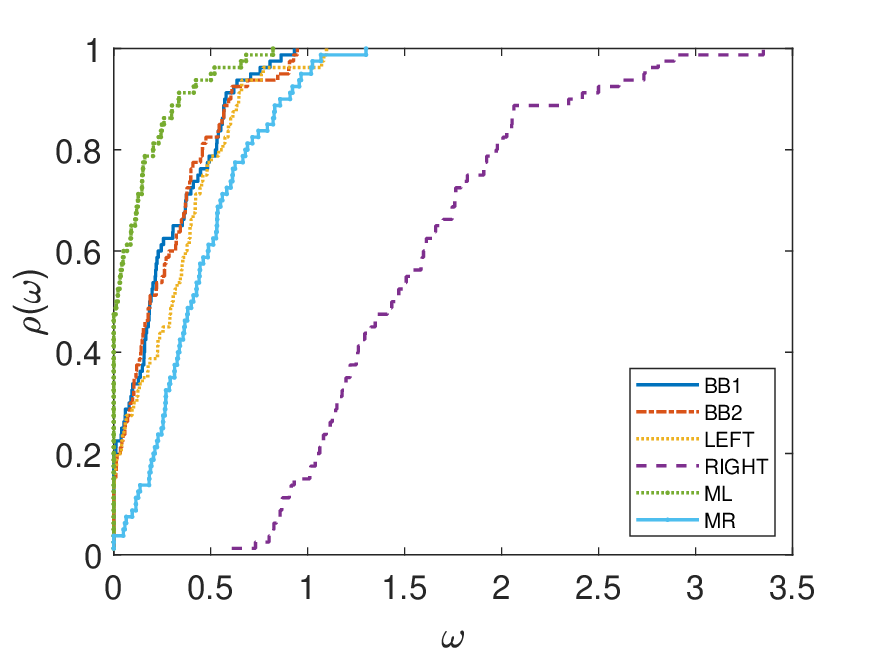}}\hspace{-14pt}
	\subfigure{
		\includegraphics[width=0.33\textwidth]{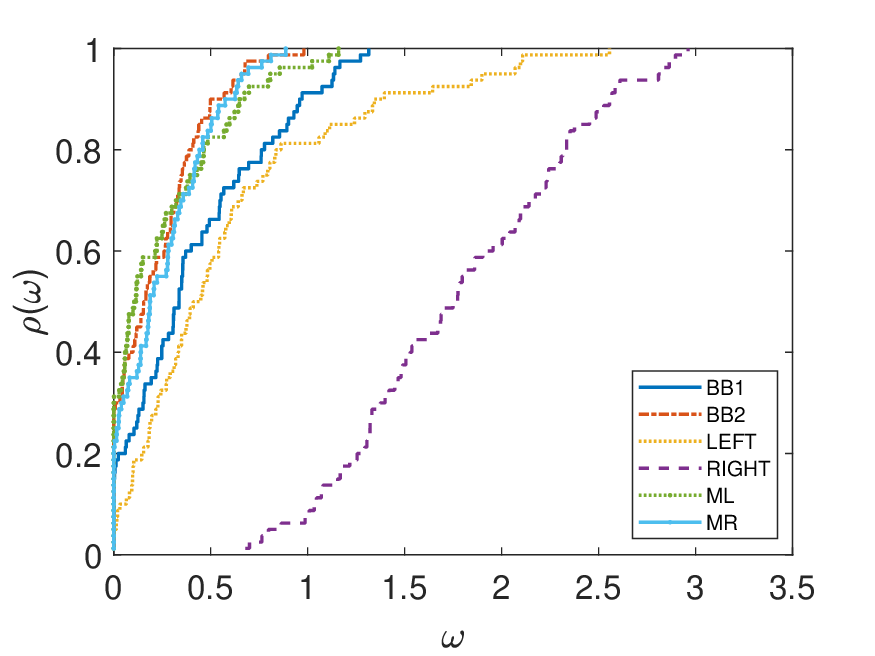}}\hspace{-14pt}
	\subfigure{
		\includegraphics[width=0.33\textwidth]{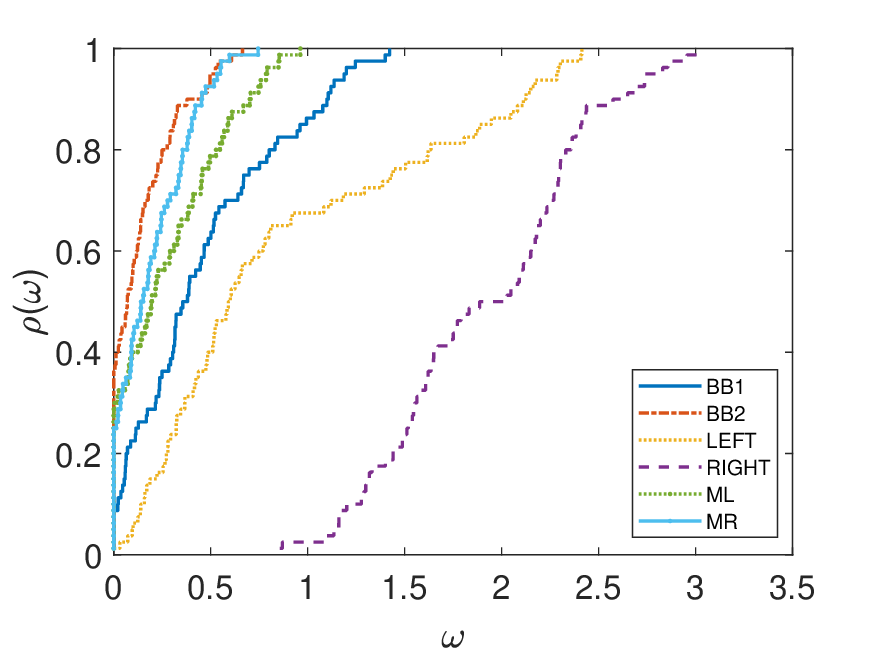}}\\\vspace{-12pt}
	\subfigure[$\varepsilon=10^{-6}$]{
		\includegraphics[width=0.33\textwidth]{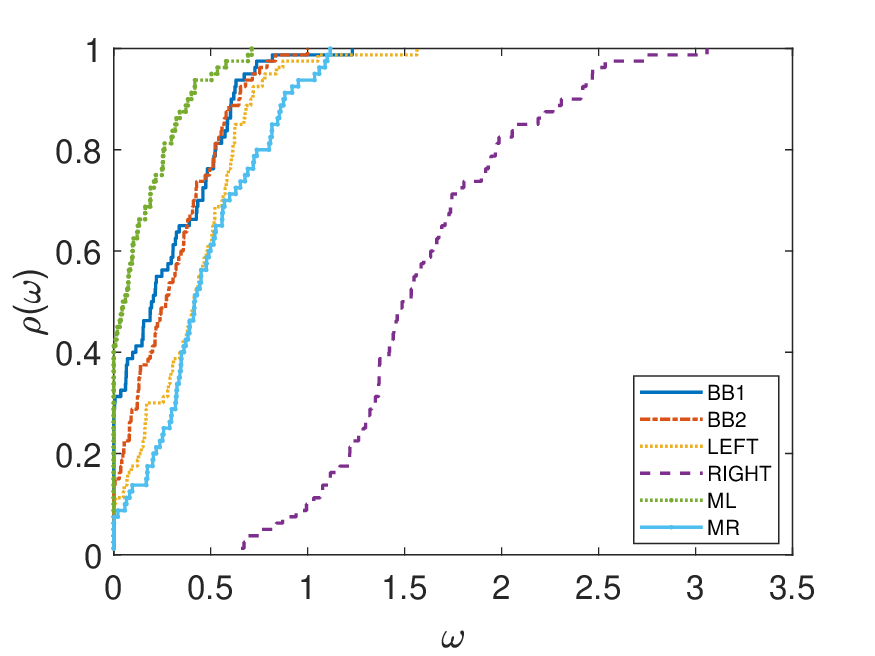}}\hspace{-14pt}
	\subfigure[$\varepsilon=10^{-9}$]{
		\includegraphics[width=0.33\textwidth]{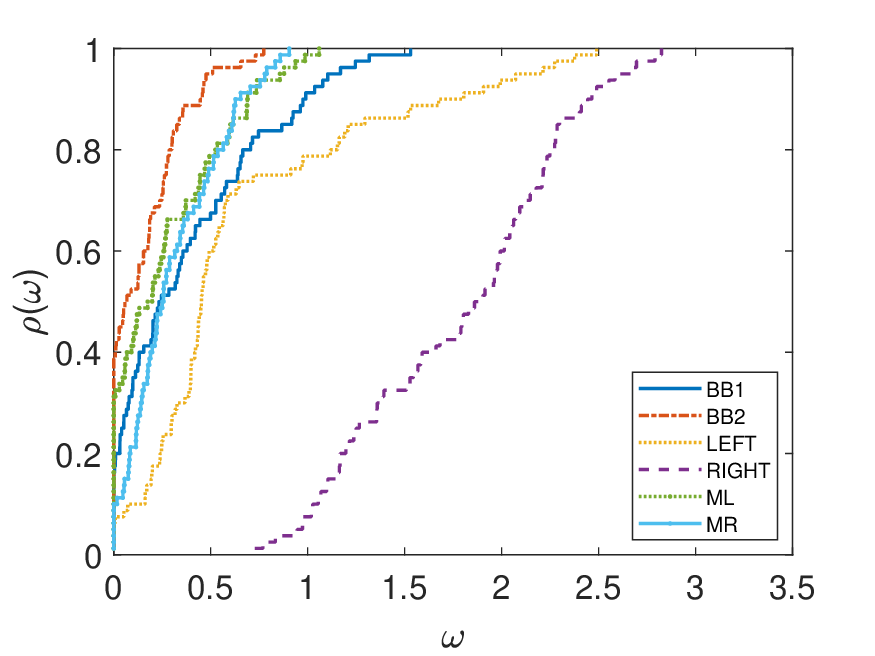}}\hspace{-14pt}
	\subfigure[$\varepsilon=10^{-12}$]{
		\includegraphics[width=0.33\textwidth]{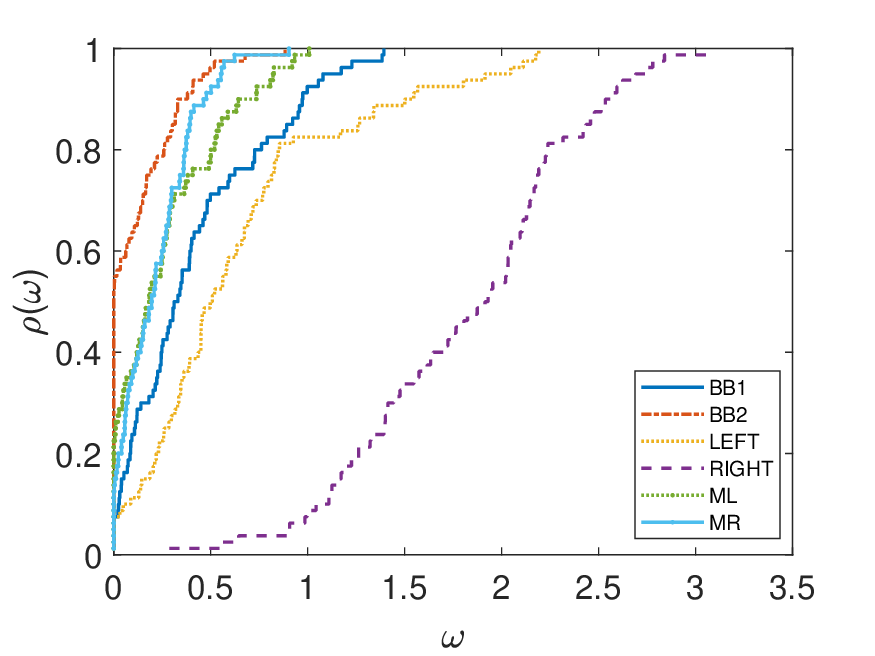}}\\
	\caption{\textit{Performance profiles of the these method on the quadratic problems \eqref{pro:qua}, from the first row to the last row, $n=500,1000,5000$, from the first column to the last column, $\varepsilon=10^{-6},10^{-9},10^{-12}$, based on iteration number.}}	
	\label{fig:nsq}
\end{figure}

We first test on the quadratic problem \eqref{pro:qua} with $A$ from \eqref{diagconditionn}. The objective function remains \eqref{pro:qua}, with problem characteristics determined by $A$. The initial point $\x_{1}=\mathbf{0}$, while the optimal point $\mathbf{x}_{*}$ has components randomly generated in $[-1,1]$. The initial step size is $\frac{\g_{1}^{\T}\g_{1}}{\g_{1}^{\T}A\g_{1}}$. Problem dimensions $n\in\{10, 50, 100, 500, 1000, 2000, 5000\}$ span small to large scales. Termination occurs when iterations exceed $20000$ or $\|\g_{k}\|_{2}<\varepsilon\|\g_{1}\|_{2}$ for $\varepsilon=10^{-6}$, $10^{-9}$, $10^{-12}$.For condition numbers $\kappa(A)\in\{10^3, 10^4, 10^5, 10^6\}$, each algorithm executes $20$ independent runs per $(\kappa, n, \varepsilon)$ combination. Figures \ref{fig:nsq1} and \ref{fig:nsq} compares algorithm performance based on iteration counts.

Figure \ref{fig:nsq1} shows that for small-scale problems ($n<100$), MR performs best. As dimension increases, MR's performance declines while ML improves. At $\varepsilon=10^{-6}$, ML consistently outperforms all other methods. For $n\geq1000$ with $\varepsilon\in\{10^{-9}, 10^{-12}\}$, BB2 performs best. These results indicate MR is sensitive to problem size, while ML exhibits greater robustness. Overall, ML outperforms BB1, but MR underperforms BB2 except for small problems.


\begin{figure}[!h]
	\centering
	\subfigure[$\varepsilon=10^{-6}$]{
		\includegraphics[width=0.33\textwidth]{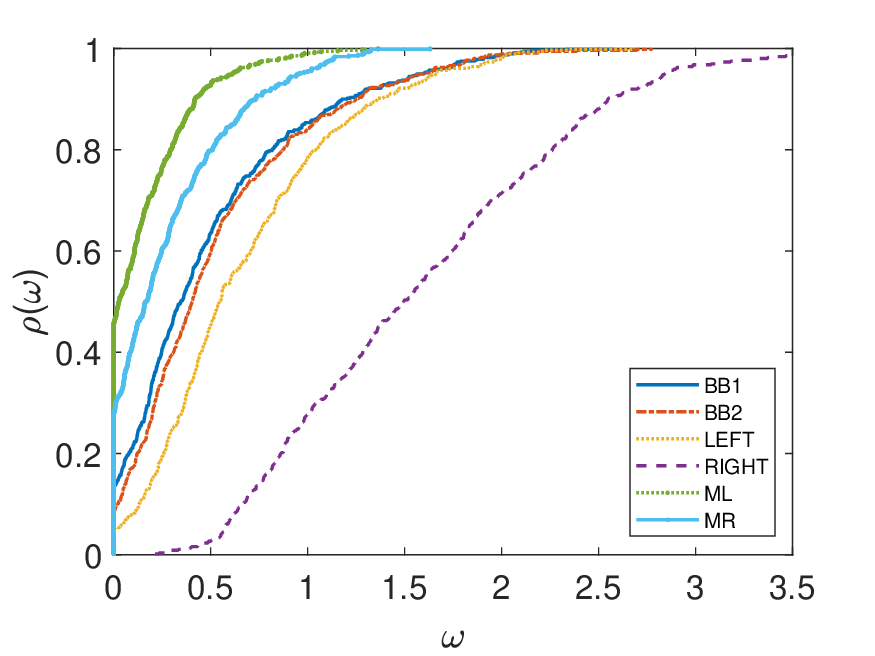}}\hspace{-14pt}
	\subfigure[$\varepsilon=10^{-9}$]{
		\includegraphics[width=0.33\textwidth]{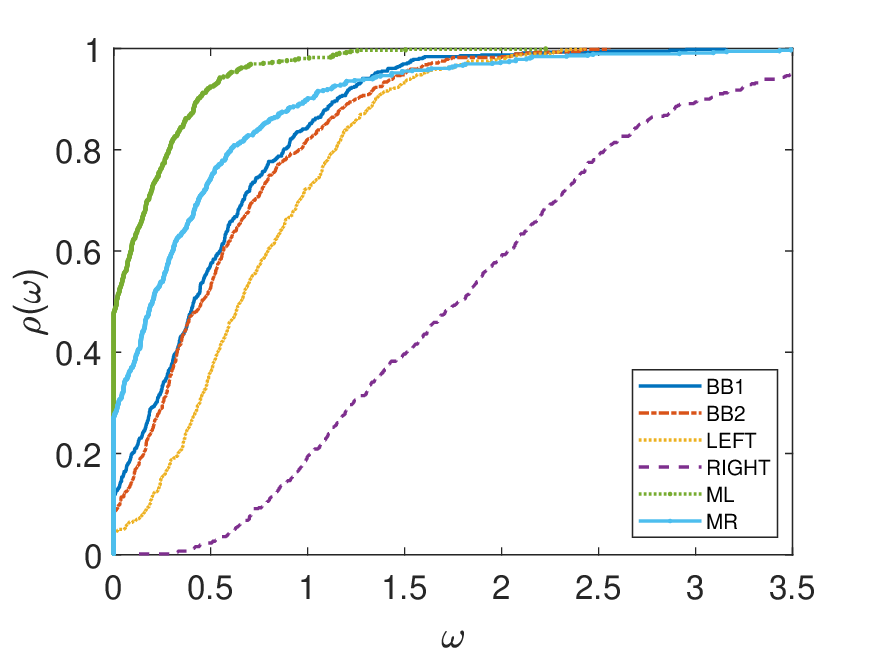}}\hspace{-14pt}
	\subfigure[$\varepsilon=10^{-12}$]{
		\includegraphics[width=0.33\textwidth]{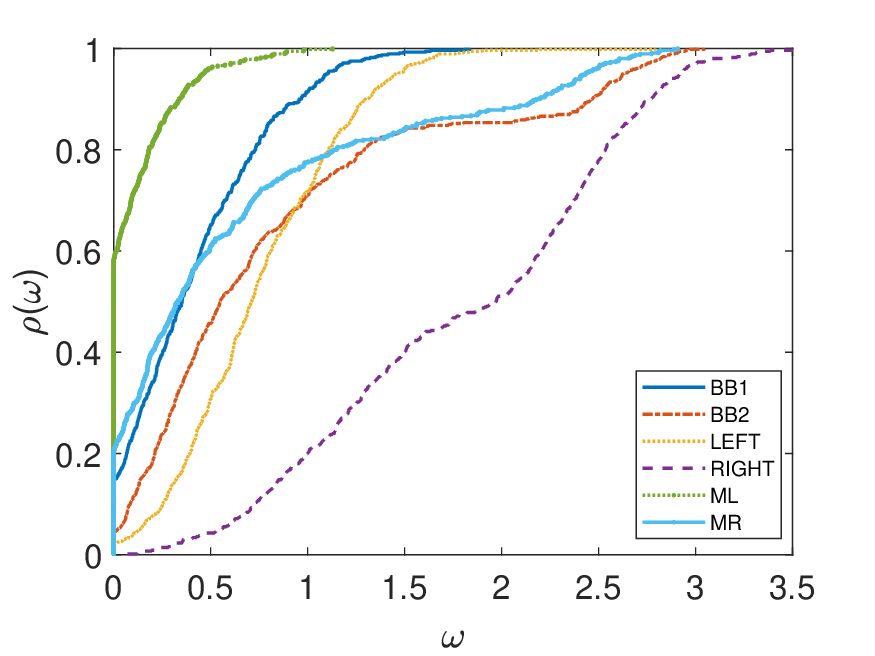}}\\
	\caption{\textit{Performance profiles of the these method on the quadratic problems \eqref{pro:qua} with different distributions of $\lambda_{j}$ from Table \ref{tab:spectrum}, $n=1000$, based on iteration number.}}	
	\label{fig:sqseven}
\end{figure}

We now consider $A$ with seven kinds of distributions of $\lambda_{j}$ summarized in Table \ref{tab:spectrum} from \cite{Dai2019familyspectralgradient}. The seven types of $A$ have different eigenvalue distributions and are often used to test the performance of spectral gradient step sizes. In these problems, we set the initial point $\x_{1}$ to be a random vector with each element in $[-5,5]$, the optimal point $\x_{*}$ is also a randomly generated vector with each element in $[-10,10]$, the problem size is $n=1000$, and the rest of the settings are the same as in the preceding test. Figure \ref{fig:sqseven} shows the performance comparison of these algorithms on these problems. 

From the results in Figure \ref{fig:sqseven}, we can see that the ML method significantly better than the other compared algorithms. When the convergence accuracy is $\varepsilon=10^{-6}$ or $\varepsilon=10^{-9}$, MR method significantly outperforms BB1 and BB2 methods. The LEFT and RIGHT methods always perform the worst. These numerical results are consistent with our analysis in Section \ref{sec:converge}.
\begin{table}[!ht]
	\centering
	\caption{Different distributions of $\lambda_{j}$ for the problem \eqref{pro:qua}}
	\setlength{\tabcolsep}{12pt}{
		\begin{tabular}{|cccc|cccc|}
			\toprule
			P     & \multicolumn{3}{c|}{$\lambda_{j}$} & P     & \multicolumn{3}{c|}{$\lambda_{j}$} \\
			\hline
			1     & \multicolumn{3}{c|}{$\{\lambda_{2},\ldots,\lambda_{n-1}\}\subset(1,\kappa)$} & \multirow{3}[4]{*}{5} & \multicolumn{3}{c|}{$\{\lambda_{2},\ldots,\lambda_{n/5}\}\subset(1,100)$} \\
			\cmidrule{1-4}    \multirow{2}[2]{*}{2} & \multicolumn{3}{c|}{$\{\lambda_{2},\ldots,\lambda_{n/5}\}\subset(1,100)$} &       & \multicolumn{3}{c|}{$\{\lambda_{n/5+1},\ldots,\lambda_{4n/5}\}\subset(100,\frac{\kappa}{2})$} \\
			& \multicolumn{3}{c|}{$\{\lambda_{n/5+1},\ldots,\lambda_{n-1}\}\subset(\frac{\kappa}{2},\kappa)$} &       & \multicolumn{3}{c|}{$\{\lambda_{4n/5+1},\ldots,\lambda_{n-1}\}\subset(\frac{\kappa}{2}, \kappa)$} \\
			\hline
			\multirow{2}[2]{*}{3} & \multicolumn{3}{c|}{$\{\lambda_{2},\ldots,\lambda_{\frac{n}{2}}\}\subset(1,100)$} & \multirow{2}[2]{*}{6} & \multicolumn{3}{c|}{$\{\lambda_{2},\ldots,\lambda_{10}\}\subset(1,100)$} \\
			& \multicolumn{3}{c|}{$\{\lambda_{n/2+1},\ldots,\lambda_{n-1}\}\subset(\frac{\kappa}{2},\kappa)$} &       & \multicolumn{3}{c|}{$\{\lambda_{11},\ldots,\lambda_{n-1}\}\subset(\frac{\kappa}{2},\kappa)$} \\
			\hline
			\multirow{2}[2]{*}{4} & \multicolumn{3}{c|}{$\{\lambda_{2},\ldots,\lambda_{4n/5}\}\subset(1,100)$} & \multirow{2}[2]{*}{7} & \multicolumn{3}{c|}{$\{\lambda_{2},\ldots,\lambda_{n-10}\}\subset(1,100)$} \\
			& \multicolumn{3}{c|}{$\{\lambda_{4n/5+1},\ldots,\lambda_{n-1}\}\subset(\frac{\kappa}{2},\kappa)$} &       & \multicolumn{3}{c|}{$\{\lambda_{n-9},\ldots,\lambda_{n-1}\}\subset(\frac{\kappa}{2},\kappa)$} \\
			\hline
	\end{tabular}}%
	\label{tab:spectrum}%
\end{table}%

\begin{figure}[!h]
	\centering
	\subfigure[$\varepsilon=10^{-6}$]{
		\includegraphics[width=0.49\textwidth]{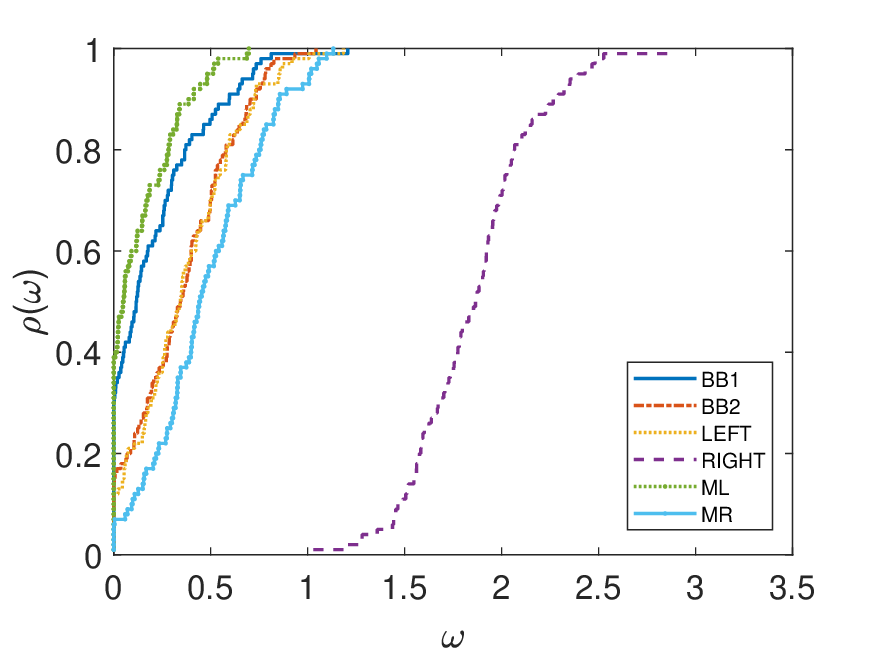}}\hspace{-18pt}
	\subfigure[$\varepsilon=10^{-9}$]{
		\includegraphics[width=0.49\textwidth]{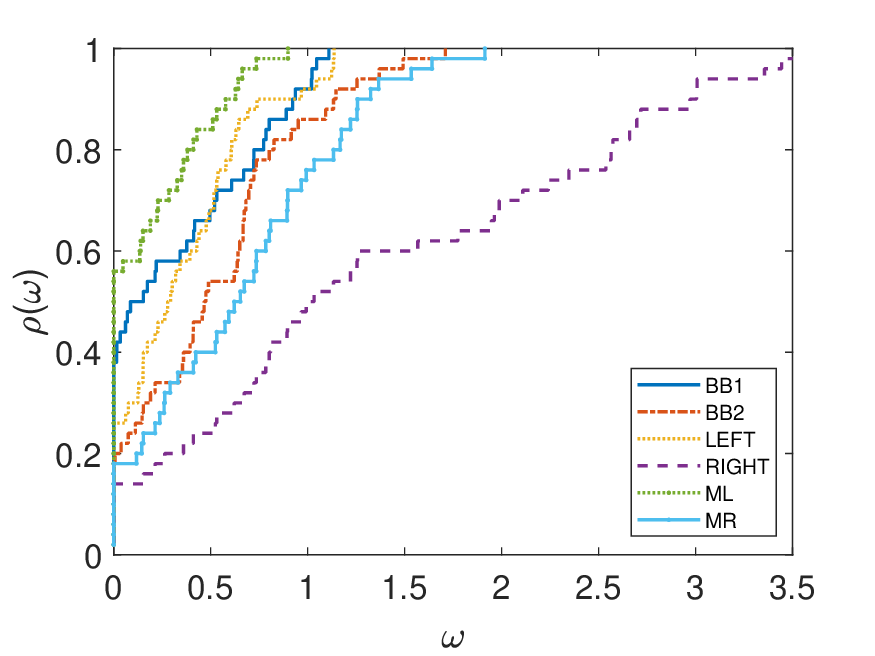}}\\
	\caption{\textit{Performance profiles of the these method on the quadratic problems \eqref{pro:qua} with $A$  \eqref{equ:TBV}, $n=1000$, based on iteration number.}}	
	\label{fig:two-point}
\end{figure}

We next compare these methods on a two-point boundary value problem \cite{Huang2022accelerationBarzilaiBorweinmethod} which can be transferred as a linear system $A\x=\mathbf{b}$ by the finite difference method. Specifically, the matrix $A=(a_{i,j})$ is given by
\begin{equation}\label{equ:TBV}
a_{i,j}=\begin{cases}
\frac{2}{h^2},\quad&\text{if}\quad i=j,\\
-\frac{1}{h^2},\quad&\text{if}\quad i=j\pm 1,\\
0,\quad&\text{otherwise},
\end{cases}
\end{equation}
where $h=11/n$. Obviously, $\kappa(A)$ increases as $n$ becomes large. Likewise, we consider the objective function \eqref{pro:qua}. The initial point $\x_{1}=\mathbf{1}$, and the optimal point $\mathbf{x}_{*}$ is randomly generated, with each component being a random number between $[-10,10]$. The scale $n$ of the problem is $1000$. In the termination conditions, we only consider $\varepsilon=10^{-6}, 10^{-9}$, and the rest of the parameter settings, including the number of independent runs, are consistent with those in the first test. Figure \ref{fig:two-point} shows the performance comparison of these methods in this two-point boundary value problem.

The results in Figure \ref{fig:two-point} show that the performance of the ML method is significantly better than that of the BB1, LEFT, BB2, MR, and RIGHT methods. 

\section{Conclusions and Discussion}
In this work, we extend the range of BB step sizes through an interpolated least-squares model, obtaining two new BB-like step sizes. We analyze the convergence and stability properties of these two methods. Intriguingly, dynamical systems analysis based on difference equations reveals that for parameter value $p\in[1.5,2)$, the LEFT method exhibits unexpected stability, characterized by monotonic decrease in the Euclidean norm of the generated gradients. In contrast, the RIGHT method proves almost invariably unstable. These results conflict with our  intuition. Intuitively, the LEFT method produces longer step sizes than BB1, which would be expected to induce more pronounced nonmonotonic behavior. Conversely, the RIGHT method yields a shorter step sizes than BB2, which should lead to a more stable monotonicity. These counterintuitive theoretical results are confirmed by numerical results. These interesting results motivate us to further explore the properties of BB step sizes.

In addition, we would like to point out that the results of the first part of the numerical experiments demonstrate that the MR method is  sensitive to the scale of the problem. The method almost presents superlinear convergence in small-scale problems, while its performance in large-scale problems is weaker than that of BB2, especially under hight-precision convergence conditions. A thought-provoking question is what causes MR to be sensitive to the scale of the problem. This question will force us to design more efficient truncation strategies in the future.

\textbf{Funding} This work was supported by the National Natural Science Foundation of China (No. 12371099).

\textbf{Data availability} Not applicable.

\section*{Compliance with Ethical Standards}

\textbf{Conflict of Interest} On behalf of all authors, the corresponding author states that there is no conflict of interest.

\noindent\textbf{Acknowledgements}
The author would like to thank Prof. Congpei An for his patient guidance.

\bibliographystyle{plain}
\bibliography{VBB}

\begin{thebibliography}{10}

\bibitem{An2025RegularizedBarzilaiBorwein}
Congpei An and Xin Xu.
\newblock {Regularized Barzilai-Borwein Method}.
\newblock {\em Numerical Algorithms}, June 2025.

\bibitem{Barzilai1988TwoPointStep}
Jonathan Barzilai and Jonathan~M. Borwein.
\newblock {Two-Point Step Size Gradient Methods}.
\newblock {\em IMA Journal of Numerical Analysis}, 8(1):141--148, 1988.

\bibitem{Bonettini2009scaledgradientprojection}
S.~Bonettini, R.~Zanella, and L.~Zanni.
\newblock A scaled gradient projection method for constrained image deblurring.
\newblock {\em Inverse Problems}, 25(1):015002, November 2009.

\bibitem{Dai2013NewAnalysisBarzilai}
Yuhong Dai.
\newblock {A New Analysis on the Barzilai-Borwein Gradient Method}.
\newblock {\em Journal of the Operations Research Society of China},
  1(2):187--198, March 2013.

\bibitem{Dai2019familyspectralgradient}
Yuhong Dai, Yakui Huang, and Xinwei Liu.
\newblock A family of spectral gradient methods for optimization.
\newblock {\em Computational Optimization and Applications}, 74(1):43--65,
  2019.

\bibitem{Dai2002Rlinearconvergence}
Yuhong Dai and Lizhi Liao.
\newblock {R-linear convergence of the Barzilai and Borwein gradient method}.
\newblock {\em IMA Journal of Numerical Analysis}, 22(1):1--10, 2002.

\bibitem{DeAsmundis2014efficientgradientmethod}
Roberta De~Asmundis, Daniela di~Serafino, William~W. Hager, Gerardo Toraldo,
  and Hongchao Zhang.
\newblock {An efficient gradient method using the Yuan steplength}.
\newblock {\em Computational Optimization and Applications}, 59(3):541--563,
  June 2014.

\bibitem{DiSerafino2018steplengthselectiongradient}
Daniela Di~Serafino, Valeria Ruggiero, Gerardo Toraldo, and Luca Zanni.
\newblock On the steplength selection in gradient methods for unconstrained
  optimization.
\newblock {\em Applied Mathematics and Computation}, 318(1):176--195, February
  2018.

\bibitem{Dolan2002Benchmarkingoptimizationsoftware}
Elizabeth~D. Dolan and Jorge~J. Moré.
\newblock Benchmarking optimization software with performance profiles.
\newblock {\em Mathematical Programming}, 91(2):201--213, 2002.

\bibitem{Elaydi2005IntroductionDifferenceEquations}
Saber Elaydi.
\newblock {\em {An Introduction to Difference Equations}}.
\newblock Undergraduate Texts in Mathematics. Springer-Verlag, New York, third
  edition, 2005.

\bibitem{Ferrandi2023harmonicframeworkstepsize}
Giulia Ferrandi, Michiel~E. Hochstenbach, and Nataša Krejić.
\newblock A harmonic framework for stepsize selection in gradient methods.
\newblock {\em Computational Optimization and Applications}, 85(1):75--106,
  February 2023.

\bibitem{Fletcher2005BarzilaiBorweinMethod}
Roger Fletcher.
\newblock {On the Barzilai-Borwein Method}.
\newblock In K.~Teo L.~Qi and X.~Yang, editors, {\em Optimization and Control
  with Applications}, pages 235--256, Boston, MA, USA, 2005. Springer US.

\bibitem{Forsythe1968asymptoticdirectionsthes}
George~E. Forsythe.
\newblock On the asymptotic directions of thes-dimensional optimum gradient
  method.
\newblock {\em Numerische Mathematik}, 11(1):57--76, January 1968.

\bibitem{Frassoldati2008Newadaptivestepsize}
Giacomo Frassoldati, Luca Zanni, and Gaetano Zanghirati.
\newblock New adaptive stepsize selections in gradient methods.
\newblock {\em Journal of Industrial and Management Optimization},
  4(2):299--312, 2008.

\bibitem{Golub2013MatrixComputations}
Gene~H. Golub and Charles~F. Van~Loan.
\newblock {\em {Matrix Computations}}.
\newblock Johns Hopkins University Press, Philadelphia, PA, fourth edition,
  2013.

\bibitem{Huang2022accelerationBarzilaiBorweinmethod}
Yakui Huang, Yuhong Dai, Xinwei Liu, and Hongchao Zhang.
\newblock {On the acceleration of the Barzilai-Borwein method}.
\newblock {\em Computational Optimization and Applications}, 81(3):717–740,
  January 2022.

\bibitem{Raydan1993BarzilaiBorweinchoice}
Marcos Raydan.
\newblock {On the Barzilai and Borwein choice of steplength for the gradient
  method}.
\newblock {\em IMA Journal of Numerical Analysis}, 13(3):321--326, 1993.

\bibitem{Xu2025ParameterizedBarzilaiBorwein}
Xin Xu.
\newblock {A Parameterized Barzilai-Borwein Method via Interpolated Least
  Squares }.
\newblock {\em arXiv}, 2025.

\bibitem{Yuan2006newstepsizesteepest}
Yaxiang Yuan.
\newblock {A new stepsize for the steepest descent method}.
\newblock {\em Journal of Computational Mathematics}, 24(2):149--156, 2006.

\bibitem{Zhou2006GradientMethodsAdaptive}
Bin Zhou, Li~Gao, and Yuhong Dai.
\newblock {Gradient Methods with Adaptive Step-Sizes}.
\newblock {\em Computational Optimization and Applications}, 35(1):69--86,
  September 2006.

\end{thebibliography}

\end{document}